\documentclass[pdf,12pt]{amsart}
\usepackage{times}
\usepackage{curves}
\usepackage{amssymb,amsfonts,amsmath,amscd}
\usepackage[pdftex]{graphicx}

\begin{document}

\newtheorem{theo}{Theorem}
\newtheorem{lem}{Lemma}
\newtheorem{cor}{Corollary}
\newtheorem{prop}{Proposition}
\newtheorem{ex}{Example}
\newtheorem{remar}{Remark}
\newtheorem{rul}{Rule}
\newtheorem{conj}{Conjecture}
\newtheorem{defi}{Definition}
\newtheorem{concl}{Conclusion}

\newcommand{\bt}{\begin{theorem}\em}
\newcommand{\et}{\end{theorem}}
\newcommand{\petet}[2]{{}^{#2}\!{#1}}
\newcommand{\side}[2]{{}^{#1}\!{#2}}

\newcommand{\mop}[2]{\stackrel{{}_{{#2}}}{#1}}
\newcommand{\Mop}[4]{\overset{\ \ \ #4}
{\underset{\ \ \ #3}
{\mathop{\sideset{^{{}_{{#2}}}}{}#1}}}}
\newcommand{\ovl}[1]{\overline{#1}}
\newcommand{\Bak}[1]{\mathbf B_{#1}}
\newcommand{\bak}{\mathbf b}
\newcommand{\dg}[2]{\mathrm{ex}_{#1}#2}
\newcommand{\ord}[2]{\mathrm{ord}_{#1}#2}
\newcommand{\Ex}[1]{\!\!\uparrow^{#1}}
\newcommand{\Exx}[2]{(#1)\!\!\uparrow^{#2}}
\newcommand{\STDm}[1]{\mathrm{STD}({#1})}
\newcommand{\STD}[1]{{\rm STD}$({#1})$}
\newcommand{\supp}[1]{{\mathrm{supp}(#1)}} 
\newcommand{\s}[1]{{#1}^{\star}}
\newcommand{\J}{\mathbf J}

\newcommand{\beas}{\begin{eqnarray*}}
\newcommand{\eeas}{\end{eqnarray*}}
\newcommand{\bea}{\begin{eqnarray}}
\newcommand{\eea}{\end{eqnarray}}
\newcommand{\lcm}{\mathrm{lcm}}
\newcommand{\md}{\mathrm{mod\ }}
\newcommand{\expon}[2]{\mathrm{ex}_{#1}\left({#2}\right)}
\newcommand{\z}{\mathbf 0}
\newcommand{\one}{\mathbf 1}
\newcommand{\rem}{\mathrm{rem}}

\newcommand{\rt}[1]{\mathrm{rem}\left(#1,\tau\right)}
\newcommand{\mulmi}[1]{\left\lfloor #1\right\rfloor_{\tau}}
\newcommand{\mulma}[1]{\left\lceil #1\right\rceil_{\tau}}
\newcommand{\mm}[1]{\,(\mathrm{mod}\ #1)}
\newcommand{\mv}[1]{\, \mathrm{mod}(#1)}

\newcommand{\pr}{{\sl Proof.\ }}
\newcommand{\bx}{\ \ \ $\Box$}
\newcommand{\bxm}{\ \ \ \Box}
\newcommand{\ind}{\mathrm{ind}\,}
\newcommand{\Acal}{\mathcal{A}}
\newcommand{\Sal}{\mathcal{S}}
\newcommand{\Eal}{\mathcal{E}}
\newcommand{\Pal}{\mathcal{P}}
\newcommand{\Hal}{\mathcal{H}}
\newcommand{\Mal}{\mathcal{M}}
\newcommand{\Ral}{\mathcal{R}}
\newcommand{\Cal}{\mathcal{C}}
\newcommand{\ccirci}[2]{\mathcal{C}(#1)\big|_{#2}}
\newcommand{\rcirci}[2]{\mathcal{R}(#1)\big|^{#2}}
\newcommand{\Bfrak}{\mathfrak B}
\newcommand{\Bbol}{\mathbf B}
\newcommand{\bbN}{\mathbb N}
\newcommand{\mb}[1]{\mathbb{#1}}
\newcommand{\m}[1]{\mathfrak{#1}}
\newcommand{\Z}{\mathbb Z}
\newcommand{\N}{\mathbb N}
\newcommand{\R}{\mathbb R}
\newcommand{\ebxt}[2]{{#1}^{\boxtimes#2}}
\newcommand{\bxt}[2]{{#1}\boxtimes{#2}}
\newcommand{\BT}{\boxtimes}
\newcommand{\setp}{\{0,\dots,p-1\}}
\newcommand{\aut}[1]{\mathbf {Aut}(#1)}
\newcommand{\debt}[2]{{#1}^{\BT{#2}}}

\newcommand{\ED}{\end{document}}

\sloppy

\title[Reversibility]{Reversibility of additive CA as function of cylinder size}

\author[V.K. Bulitko]{Valeriy Bulitko\\
\ \\
{\tiny Athabasca University}\\
{\tiny valeriyb@athabascau.ca}}

\vspace*{-2.5cm} \sloppy \thispagestyle{empty}

\maketitle

\begin{abstract}
Additive CA on a cylinder of size $n$ can be represented by 01-string $V$ of length $n$ which is its rule. We study a problem: a class $S$ of rules given, for any $V\in S$ describe all sizes $n', n'>n,$ of cylinders such that extension of $V$ by zeros to length $n'$ represents reversible additive CA on a cylinder of size $n'$. Since all extensions of $V$ have the same collection of positions of units, it is convenient to say about classes of collections of positions instead of classes of rules.  A criterion of reversibility is proven. The problem is completely solved for infinite class of ''block collections'', i.e. $\{(0,1,\dots,h)|h\in\Z^+\}$. Results obtained for  ``exponential collections'' $\{(1,2,4,\dots,2^h)|h\in\Z^+\}$  essentially reduce the complexity of the problem for this class. Ways to transfer the results on other classes of rules/collections are described. A conjecture is formulated for class $\{(0,1,2^m)|m\in\Z^+\}$.\\

\noindent {\bf MSC-class:} 37B15, 68Q80

\end{abstract}

The problem of reversibility of general cellular automata comes from a general theory of computations as well as from applications.  In algorithmic formulation the problem was shown to be undecidable for dimensions higher 1 and decidable in dimension 1. There are more subtle aspects of those results relating to complexity (descriptive and algorithmic). See an overview and some references in\\
\verb= http://en.wikipedia.org/wiki/Reversible\_cellular\_automaton=

On the other hand consideration of semi-groups of automata also leads to reversibility problems of a certain kind. In particular, if we deal with additive CA on a cylinder of size $n$ it is possible to represent automata rules and states by elements of a semi-group built on strings of length $n$ with a discrete convolution as an operation. From this prospect reversible automata create a subgroup of the semi-group. If we want to general features of the groups of reversible additive CA on cylinders we might want to find a way for any size $n$ of a cylinder to construct reversible strings that could constitute a set of generators  of the group for the cylinder. 

Let $\mathfrak M_n$ be the set of all 01-strings of length $n$. For any $V,W\in\mathfrak M_n$ we denote $j$-th component of $V$ by $V(j),0\le j<n,$ and denote convolution of $V,W$ as $V\BT W$ where $[V\BT W](j)=\sum_{s=0}^{n-1}V(s)W(j-s)\ \mm 2$. $\mathfrak S_n$ is a semi-group on $\mathfrak M_n$ created by operation $\BT$. Let $\mathfrak G_n$ be a sub-group of $\mathfrak S_n$ of all reversible elements. 

In a slightly more general form the problem might be formulated as following: a subset $S$ of $\Z^+$ and a function $f:S\to \mathfrak M_n$ given, describe $\{n\in S| f(n)\in \mathfrak G_n\}$. 
Although in this paper we consider probably the simplest couples $(S, f)$. Namely let $V\in\mathfrak M_n$. We define $S_V,f_V$ as following: $S_V=\{n'\in\Z^+|n'>n\},\  f_V(n') $ is a string that has length $n'$ and s.t. $V$ is its prefix whereas all positions $i, n\le i<n',$ are occupied with 0. In other words we consider all extensions of $V$ with zeros and search for which lengths  the extensions are reversible. For example, if $V=(1,1,0,0,1)$ the question is:  for which $n, n>5,$ the string $(1,1,0,0,1,0,\dots,0)$ of length $n$ is reversible on a cylinder of size $n$.

A string $W$ given its $j$-th position is denoted as $W(j)$. In case when a string is represented by a complex expression like $V\BT W$ for instance, $j$-th position will be denoted as $[V\BT W](j)$. 
Expression $[W]_a^b$ is applicable to string $W$ of length $n$ when $0\le a\le b\le n$ and denotes the sub-string of $W$ that begins with position $a$ and ends with position $b$. So, $[W]_j^j=W(j)$. If denotation $W$ is a simple we sometime use $W\mid_a^b$ instead of $[W]_a^b$. 
Sometimes it is convenient to write a string $(i_1,\dots, i_k),i_j\in A,$ as a word in alphabet $A$. For example, string $(0,1,1,1,0,\dots,0)$ of length 15 could be represented by word $01^30^{11}$. 
$\sigma^j$ denotes a cyclic shift to right; namely 
\beas
[\sigma^j(V)](i)=V((i-j)\mm n).
\eeas
For any length $\z$ denotes a zero string $(0,\dots,0)$. 

Relative positions of units in a string and a cylinder size determine reversibility of the string. Let $V[x_1,\dots,x_r;\,n]$ is a string of length $n$ with $r,r\ge1,$ units such that $x_i, i= 1,\dots,r,$ are positions of units in the string where $0\le x_1<x_2<\dots<x_r\le n$.

 Given a string $V[x_1,\dots,x_r;\,n]$ let $\partial V$ denote a derivative string $V[x_2-x_1,\dots,x_r-x_1;n]$. In case $r\le1$ we assume $\partial V=\z$.
\begin{lem}\label{l0}
$V[x_1,\dots,x_r;\,n]$ given,  the following equivalences hold
\beas
V\BT L =\z \iff (\partial V)\BT L = L \iff L =\sum_{i=2}^r \sigma^{x_i-x_1}(L) \mm 2,
\eeas
where  $L$ is non-zero string of length $n$.
\end{lem}
\pr 
First
\bea\label{lo1}
[V\BT L](j)=\sum_{s=1}^nV(s)L(j-s)=\underset{s: V(s)=1}{\sum}L(j-s)=\underset{i}{\sum}L(j-x_i) = 0.
\eea
The latter equation can be represented as 
\beas
L(j-x_1)+\sum_{i=2}^rL((j-x_1)-(x_i-x_1))=0, j=0,\dots,n-1.
\eeas
Since $x_1$ is fixed, by modulo $n$ the amount $j-x_1 $ runs over the whole set $\{0,\dots,n-1\}$. Therefore the latter equation could be rewritten as 
\bea\label{lo2}
L(j)=\sum_{i=2}^rL(j-(x_i-x_1)), j=0,\dots,n-1.
\eea
Now, as $L(j-a)=[\sigma^a(L)](j)$, the latter system of equations could be rewritten as
\beas
L =\sum_{i=2}^r \sigma^{x_i-x_1}(L) \mm 2.
\eeas
On the other hand, 
\beas
[(\partial V)\BT L](j) = \sum_{s}[\partial V](s)L(j-s)
\eeas
and taking into account that $[\partial V](s)=1$ iff  $s=x_i-x_1$ for some $i\in\{2,\dots,r\}$ we get that for all $j$ it holds
\beas
[(\partial V)\BT L](j) = \sum_{s}[\partial V](s)L(j-s)= \sum_{i=2}^rL(j-(x_i-x_1)) \mm 2.
\eeas
From (\ref{lo2}) it follows
\beas
(\partial V)\BT L=L.
\eeas
\bx

\begin{lem}\label{l1}
 Let 
$x_l\in\{0,1,\dots,n-1\},l=\ovl{1,r},\ l<l'\implies x_i<x_{l'}$. Then $V[x_1,\dots,x_r;n]$ is irreversible on a cylinder of size $n$ iff a non-zero string $L$ of length $n$ exists satisfying the condition
\beas
L(j)=\sum_{i=2}^rL(j-(x_i-x_1)), j=0,\dots,n-1.
\eeas 
\end{lem}
\pr As $V\BT(T+L)=(V\BT T) + (V\BT L)$\footnote{All sums are done by modulo 2.}, it holds that 
$V[x_1,\dots,x_r;n]$ is irreversible on a cylinder of size $n$ iff there exists a non-zero string $L$ of the length $n$ such that $V\BT L=\z$. As it is shown by (\ref{lo1}), (\ref{lo2}) in the proof above 
\beas
[V\BT L](j) = 0 \iff L(j)=\sum_{i=2}^rL(j-(x_i-x_1)) \mm2
\eeas
for all $j=0,\dots,n-1$. 

We could start with another form of irreversibility of $V$. It is an existence of a linear combination of  $\sigma^i(V)$ (where $\sigma^j$ is shift to right on $j$ positions by modulo $n$) of $V[x_1,\dots,x_r;n]$ that is equal to zero. We can represent the combination by a non-zero vector $L$ of length $n$ where units show the lines of a circulant $\Cal(V)$ for  $V[x_1,\dots,x_r;n]$ occurring in the combination:
\beas
\sum_{L(j)\neq0}\sigma^j(V[x_1,\dots,x_r;n])=(0,\dots,0) \mm2.
\eeas
Let 
\beas
[K[x]](i)=
\begin{cases}
1,& i=x,\\
0,& i\neq x.
\end{cases}
\eeas
With this denotation it is possible to rewrite the LHS as
\beas
\sum_{L(j)\neq0}\sigma^j\left(\sum_{l=1}^rK[x_l]\right)=
\sum_{L(j)\neq0}\sum_{l=1}^r\sigma^j(K[x_l])= 
\sum_{l=1}^r\left[\sum_{L(j)\neq0}\sigma^j(K[x_l])\right]
\eeas
On the other hand $\sum_{L(j)\neq0}\sigma^j(K[x])=\sigma^x(L)$. Indeed, $\left[\sum_{L(j)\neq0}\sigma^j(K[x])\right](i)=1\iff\exists j[L(j)=1\ \&\  x+j=i]\iff [\sigma^x(L)](i)=1$. Therefore 
\beas
\sum_{L(j)\neq0}\sigma^j\left(\sum_{l=1}^rK[x_l]\right)= \sum_{l=1}^r\sigma^{x_l}(L)
\eeas 
Thus 
\beas
\sum_{L(j)\neq0}\sigma^j(V[x_1,\dots,x_r;n])=(0,\dots,0) \mm 2 \iff \sum_{l=1}^r\sigma^{x_l}(L)=[0,\dots,0] \mm2.
\eeas
Hence for each $j\in\{0,\dots,n-1\}$ we have $\sum_{l=1}^r\sigma^{x_l}(L)|_j^j=0\mm2$. Since $\sigma^a(T)|_b^b=T(b-a)$, we arrive at 
\beas
\sum_{l=1}^rL(j-x_l)=0 \mm2,\ j=0,\dots,n-1.
\eeas
This is the last equality from (\ref{lo1}). \ \ \  \bx

\begin{lem}\label{l2}
If $r,x_1,\dots,x_r$ are arbitrary non-negative integers s.t.  $r>1, x_1<x_2<\dots<x_r$ and $\omega$ is any $(0,1)$-string of length $\delta,\delta=x_r-x_1$, then there exists a unique unlimited from the right end $(0,1)$-string $T$  satisfying the recursion
\bea\label{RR1}
\begin{cases}
\forall i<\delta & [T(i)=\omega(i)],\\
\forall i \ge \delta & [T(i)=\sum_{l=2}^rT(i-(x_l-x_1))].
\end{cases}
\eea 
The solution $T$ is periodic after a pre-period. For the lengths $\lambda_{pp},\lambda_{p}$ of the pre-period and period respectively,  $\lambda_{pp}+\lambda_p\le |2|^{\delta}-1$. 
\end{lem}
\pr 
$\omega$ given, $T$ is defined by the recursion uniquely. 
If $\omega$ is zero-string, obviously $T$ is also zero-string with period of length 1 and empty pre-period. 
Let $\omega$ be non-zero. From the recursion~\ref{RR1} it follows that apart from the prefix $\omega$ each sub-string (sub-word)  of $T$ of length $\delta$ depends only on the previous sub-string of length $\delta$. When sub-strings $T|_t^{t+\delta-1}$ and $T|_{t+k}^{t+k+\delta-1}$ coincide, the sub-string $T|_t^{t+k-1}$ appears to be a period of string $T|_t^{\infty}$. That is the solution should be periodical starting from some point $t$.  Since there could be no more than $|2|^{\delta}$ the sum of lengths of the pre-period and period do not exceed $|2|^{\delta}-1$. \bx

\begin{cor}\label{cor3}
The recursion~(\ref{RR1}) has a solution $T$ for some string $\omega$ iff (\ref{RR1}) has purely periodic solution $T'$  with the same period.
\end{cor}
\pr
Indeed, if $T$ is a solution unlimited from right with a prefix $\omega$ then for any integer $t$ the string $T'=T|_t$ also is a solution to the recursion but maybe with another prefix $\omega'$ where $\omega'=T_t^{t+\delta}$. 
\bx

Examples $(x_1,x_2,x_3)\in \{ (1,2,4),(1,3,4),(1,4,5),(1,2,5),(1,4,8),(1,3,12),(1,2,16),$ $(1,8,16)\}$ (see Table~\ref{tab1} below, lines 3,4,8,10,11,12,15,21) show that  the upper bound $|2|^{\delta}-1$ for $\lambda_{p}$ is reachable. 

An odd $r$ and positions $x_1,\dots,x_r$ given, the point of interest is the spectrum of all lengths of  periods of non-zero solutions to~(\ref{RR1}) for possible $\omega$ because the knowledge allows to conclude about reversibility of $V[x_1,\dots,x_r;n]$. From this angle of view it makes a sense to introduce the following terms. We call {\sl collection of positions} any not empty finite sequence (in increasing order) of non-negative integers. If $C=(x_1,\dots,x_r)$ is a collection of positions we call collection of {\sl shifts} for $C$ the derivative collection $\partial C$  that is $(r-1)$-tuple $(x_2-x_1,\dots,x_r-x_1)$. So, elements of $\partial C$ are shifts. Note that the recursion~(\ref{RR1}) is defined in terms of the derivative collections, i.e. in terms of shifts.

\begin{theo}\label{t2}
 A collection $C=(x_1,\dots,x_r)$ of positions given:
\begin{itemize}
\item [(i)]  the string 
$V[x_1,\dots,x_r;n]$ is irreversible\footnote{Naturally on a cylinder of size $n$.} iff $n$ is multiple of the length of a period of a possible solution to the recursion~(\ref{RR1}); 
\item [(ii)] If $V[x_1,\dots,x_r;n]$ is irreversible, $n$ is multiple of a number not exceeding $2^{x_r-x_1}-1$; 
\item [(iii)] when $r=1$ the string $V[x_1;n]$ is reversible for any $n>x_1$; 
\item [(iv)] if $r$ is even, $V[x_1,\dots,x_r;n]$ is irreversible for any $n,n>x_r$. 
\end{itemize}
\end{theo}
\pr If $V[x_1,\dots,x_r;n]$ is irreversible on a cylinder of size $n$ then by lemma~\ref{l1} there exist a vector $L$ of length $n$ satisfying the condition 
\bea\label{eq1-theo2}
L(i)=L(i-(x_2-x_1))+\dots+L(i-(x_r-x_1)), i=0,\dots,n-1, \mm2.
\eea 
This is the second relation of the recursion~(\ref{RR1}). To satisfy the first we set $\omega=L|_0^{x_r-x_1}$. In this case the pre-period is empty and the length of the period divides $n$.

Conversely, if recursion~(\ref{RR1}) is true for a string $L$ and its period divides $n$ we can replace $L$ with purely periodic solution $L'$ having the same period length, corollary~\ref{cor3}. As for $L'$ the relation 
\beas
L'(i)=L'(i-(x_2-x_1))+\dots+L'(i-(x_r-x_1)), i=0,\dots,n-1, \mm2,
\eeas 
is satisfied, lemma~\ref{l1} is applicable. 

The case $r=1$ is obvious because by a cyclic shift $V[x_1;n]$ is transformable into $\one$ for cylinder of size $n$. When $r$ is even one might set $L=\one$ on the cylinder of size $n$. Note, that if $r=0$ then $V=\z$.
\bx

\footnotesize
\begin{center}
 \begin{table}[h]
\caption{Some examples of collections of positions and lengths of the corresponding proper  periods.}\label{tab1}
\begin{tabular}{|r||l|l||r||l|l||} 
\hline
\# & collections & lengths of proper  periods&\#& collections & lengths of proper periods \\
\hline\hline
 1 & 1,\,2,\,3 & 3 & 13& 1,\,3,\,10& 365,\,31,\,15\\
 2& 2,\,3,\,4& 3 &  14& 1,\,3,\,11& 42,\,14,\,21,\,7,\,6,\,3\\
 3 & 1,\,3,\,4& 7 & 15& 1,\,3,\,12& 2047\\
4 &  1,\,2,\, 4& 7  & 16& 1,\,3,\,13& 126,\,63\\
 5&  3,\,4,\,5&  3 &  17& 1,\,3,\,14& 1785,\,255,\,21,\,7,\,3\\
 6&  2,\,4,\, 5&  7 &18& 1,\,3,\,15& 254,\,127\\
 7&  2,\,3,\,5&  7 & 19& 1,\,3,\,16& 4599,\,511,\,63\\
8&   1,\,4,\,5 & 15 & 20& 1,\,4,\,16& 63,\,21,\,9,\,7\\
9&   1,\,3,\, 5& 6,\,3 & 21& 1,\,8,\,16& 32767\\
 10& 1,\,2,\,5& 15 & 22& 2,\,5,\,7,\,8,\,9& 42,\,21,\,7,\,6,\,3\\
 11& 1,\,4,\,8& 127 & 23& 2,\,4,\,6,\,7,\,9& 105,\,15,\,7\\
 12& 1,\,2,\,16& 32767 & 24& 1,\,6,\,7,\,8,\,9& 217,31,7\\
\hline
\end{tabular}
\end{table}
\end{center}
\normalsize

A 01-string with finite number of units defines a collection $C$ of positions $(x_1,\dots,x_r)$. Its derivative collection of shifts determines a series $S$ of periods which we call {\sl spectrum} of both C and $\partial C$.  As we saw, any spectrum consists of numbers multiple of a finite set of numbers that are lesser  $(x_r-x_1)(2^{x_r-x_1}+1)$ (see theorem~\ref{t2}, (ii)). 
On the other hand, in fact only a maximal sub-system $K(S)$ elements of a spectrum $S$ s.t. $\forall x\in S\forall y\in K(S)[x\neq y\implies x\nmid y]$ is essential. For example, for $S=\{k\pi|\pi\in\{42,21,14,7,3\},k\in\Z^+\}$ the sub-system is $[7,3]$. Let us call $K(S)$ the {\sl spectrum kernel}. The definition $K(S)=\{y\in S| \forall x\in S[x\neq y\to x\nmid y]\}$ is correct since define the set uniquely. Namely $K(S)$ is the set of all minimal elements in partial order $x\preceq y$ defined as $x\mid y$ on the set of period lengths. 

We call {\sl block-collection} a collection of positions of kind $(a,a+1,\dots,a+h)$ where $a,h$ are some non-negative integers. 

\begin{theo}\label{blocks-t}
For any block-collection $C=(0,1,\dots,h),h\in\bbN,h>0,$ the lengths of proper periods are exactly all divisors of $h+1$ apart from 1 if $h+1$ is odd and apart from periods of length $2$ if $2|h\wedge 4\nmid h$. In case of period length 1 when $h$ is even we necessarily deal with zero-solution.  
\end{theo} 
\pr
Let  $L$ be a string satisfying the recursion with the collection $(0,1,\dots,h)$ of positions. Let us for sake of convenience the enumeration of positions of symbols in $L$ start with 1. 
We show now that $L(h+1)=\rho(L_1^h)$ (where $\rho(X)$ is the parity of a string $X$) and $L|_{h+2}^{2h+1}=L_1^h$. 

Indeed, $L(h+1)=\sum_{k=1}^{h}L(k)$ due to the fact that $\partial C=(1,\dots,h)$. And obviously $\rho(L_1^h)=\sum_{k=1}^{h}L(k)$.\footnote{Reminder: all sums should be calculated by modulo 2.} 
From here for $L(h+2)$ we have equation 
\beas
L(h+2)=L(h+1)+\sum_{k=2}^{h}L(k)=\\L(h+1)+\rho(L|_1^h)+L(1)=\rho(L|_1^h)+\rho(L|_1^h)+L(1)=L(1).
\eeas
Hence $L(h+2)=L(1)$. Assume it has been already deduced that 
\bea\label{t3-1}
1<j<h\implies L(h+j)=L(j-1). 
\eea
 From here again by direct observation of the given recursion relation for $L(h+j+1)$  we obtain:
\beas
L(h+(j+1))=\sum_{k=j+1}^{h}L(k)+\rho(L_1^h)+\sum_{k=1}^{j-1}L(k)=\left(\sum_{k=1}^{j-1}L(k)+\sum_{k=j+1}^{h}L(k)\right)+ \rho(L_1^h)=\\
\left(\rho(L_1^h)+L(j)\right)+\rho(L_1^h)=L(j).
\eeas
For $L(2h+1)$ we have recursion equation 
\beas
L(2h+1)=\sum_{k=1}^nL(h+k)=\rho(L_1^h)+\sum_{k=2}^hL(h+k)=\rho(L_1^h)+\sum_{k=2}^hL(k-1)=\\
\rho(L_1^h)+\sum_{k=1}^{h-1}L(k)=\rho(L_1^h)+\left(\rho(L_1^h)+L(h)\right)=L(h).
\eeas
Thus the sequence $L(1),\dots,L(h),\rho(L_1^h)$ repeats for any $L|_1^h$. Hence, no proper period exists with length exceeding $h+1$.

Now, let $h>1$. Obviously we can choose $L|_1^h$ so that $L|_1^h,\rho(L_1^h)$ became non-periodic and thereby realize a period of length $h+1$. For that (for instance) set $L(i)=0, 1\le i<h, L(h)=1$. In this case $\rho(L_1^h)=1$ and period looks like $0\dots011$.\footnote{Of course there could be many other settings.}
If $h=1$ the only periods are $0$ and $1$ of length 1. 

Let $h>1$ and is even. For any proper positive divisor $m$ of $h+1$ (i.e. $m< h+1$)  we can construct a period $\pi$ of length $q$ where $q=\frac{h+1}m$. Indeed, $q\ge3$ and is odd (as well as $m$ because of $h+1$ is odd). And we define $\pi=1^{q-1}0$. Since $q-1$ is even, $\sum_{k=1}^{h}L(k)=0$ when $L|_1^{h+1}$ is compiled from $\pi$ as $\pi^m$. In this construction $L(h+1)$ appears to be the last symbol of the last occurrence of $\pi$ in $L|_1^{h+1}$ and therefore is equal to $\rho(L|_1^h)$ as necessary (see above). Since $\pi$ is aperiodic, no lesser period exists for this definition of $L$. 

If $m=h+1$ we have $q=1$. The only possible solution for $L$ must consist of zeros: otherwise $L$ can consist only of ones but this is inconsistent with $\rho(L_1^h)=0$ (because $h$ is even).

Finally, let $h>1$ and is odd. In this case obviously there exist both trivial solutions when $L$ is built only from zeros or only from ones. The period is 1 in these cases. When $m $ is a proper divisor of $h+1$ (i.e. $1<q<h+1$) we have 
the above construction applicable as well if $q>2$. It is also applicable when $q=2\wedge 4\mid h+1$. Indeed, we have even $m$ for this case we set $L|_1^h=(10)^{m-1}1$ getting $\rho(L|_1^h)=0$ due to the fact that $(m-1)$ is odd. Hence $L|_1^{h+1}=(10)^m$ and 01 is a period.

For the case $q=2\wedge 4\nmid h+1$ no period of length 2 is possible: in this case $q=2, m$ is odd and we have two opportunities to try. The first $L|_1^h=(10)^{m-1}1$ yields $\rho(L|_1^h)=1$ because $m-1$ is odd and period $10$ is destructed by $L(h+1)=1$. The second  $L|_1^h=(01)^{m-1}0$ yields $\rho(L|_1^h)=0$, i.e. a sub-string 00 in $L$ that contradicts to the form of the period. \footnote{We can directly set $L(0)=0$ and $L(i)=1,i=\overline{1,2h}$ for the maximal length of the period. Another definition for the maximal period is $L(2h-1)=L(2h)=1, L(i)=0,i=\overline{0,2h-2}$.}
\bx 

 Let $C$ be a collection $x_1,x_2,\dots,x_r$ of positions, $m$ a positive integer, and $W$ be a string of length $g\ge \delta$ where $\delta=x_r-x_1$. We denote $C[m]\{W\}$ a string $T$ of length $\delta+m$ obtained by recursion~\ref{RR1}
\bea\label{RR2}
\forall i\in \{0,\dots,g+m\}
\begin{cases}
0\le i< g \implies T(i)=W(i),\\
 g\le i \implies T(i)=\sum_{l=2}^rT(i-(x_l-x_1)).
\end{cases}
\eea 
We write $C\{W\}$ in case $m=1$. And $C^a\{W\}=\underbrace{C\{C\{\dots C\{W\}\dots\}\}}_{a\text{ times         }}$.
Clearly $C[m]\{W\}=C^m\{W\}$. We define also that $C[0]\{W\}=C^0\{W\}=W$. 

\begin{theo}\label{l-sup} 
$  C[m]\{W+V\}=  C[m]\{W\}+  C[m]\{V\}$.
\end{theo}

\pr 
First we have $[ C[m]\{W+V\}]|_0^{g-1}=[W+V]|_0^{g-1}=W|_0^{g-1}+V|_0^{g-1}=[C[m]\{W\}]|_0^{g-1}+[ C[m]\{V\}]|_0^{g-1}$. Assume that $g\le i\le m+g$ and for for all $j<i$ the equality 
\beas
[C[m]\{W+V\}](j)=[C[m]\{W\}](j)+[C[m]\{V\}](j)
\eeas
is proven. Let $T= C[m]\{W\}$ and $H= C[m]\{V\}$. From~(\ref{RR2}) we have 
\beas
[ C[m]\{W+V\}](i)=\sum_{l=2}^r [C[m]\{W+V\}](i-(x_l-x_1))=_{\text{by induction hypothesis }}\\
\sum_{l=2}^r[ C[m]\{W\}+ C[m]\{V\}](i-(x_l-x_1))=\\
\sum_{l=2}^r [C[m]\{W\}](i-(x_l-x_1)) + \sum_{l=2}^r [C[m]\{V\}](i-(x_l-x_1))=T(i)+H(i).
\eeas
\bx
\begin{prop}\label{pr-sum}
Let $T,H$ are periodic strings on the same segment $I$ with lengths of periods $t,h$ respectively and the length of $I$ is multiple of $\lcm(t,h)$. Then $T+H$ has a period on $I$ whose length divides $\lcm(t,h)$. 
\end{prop}  
\pr Obviously $\lcm(t,h)$ is a period for $T+H$.  \bx 
\begin{prop}\label{sum powers}
Let $x,y,x\in\Z^+\cup\{0\}$. The equation $2^x=2^y+2^z$  has unique solutions $x=z+1,y=z$.
\end{prop}
\pr  Sure, $x=y+z,y=z$ is a solution to $2^x=2^y+2^z$. Assume $y>z$. In this case $2^x>2^y$ and therefore $2^{x-1}\ge 2^y\ \wedge\ 2^{x-1}>2^z $. Hence $2^x=2^{x-1}+2^{x-1}>2^y+2^z$. \bx

Let $\mathcal E_n=(1,2,4,\dots,2^n)$, i.e.  $i$-th shift is $s_i=2^i,i=\overline{1,n}$, and $K_{m,n}=0^{m-1}10^{n-m}$.
It is worth to note that the action of operator $C[m]$ on a string $W$ as it is defined depends only on $\partial C$ in fact. 
Therefore $\Eal_n$ as operator is defined completely by shift collection $\partial\Eal_n$ whose elements we denote below as $s_i,s_i=2^i-1, i=1,\dots,n$. We call collections $\Eal_n$ {\sl exponential collections} of positions. 

Theorem~\ref{l-sup} reduces behaviour of operators on strings to their behaviour on constituents $K_{m,n}$.

\begin{lem}\label{l7}
The following holds\footnote{For convenience we assume supports of strings to start with number 1.}
\bea\label{form-lem7}
\begin{cases}
\Eal_n^{2^{n+1}-1}\{K_{m,2^n-1}\}=\\
K_{m,2^n-1}\cdot K_{\eta,\eta}\cdot \Pi_{j=k}^{n } 
 \left[\left(\Pi_{i=0}^{k-2}K_{1,2^i}\right)\cdot K_{1,2^j-2^{k-1}+1-m\,\delta_{n,j}}\right]\cdot K_{m,2^{n}-1}
 \end{cases}
\eea
 where ``$\cdot$'' means concatenation of words, $k=\mu z[s_z\ge 2^n-m], \eta=m - (2^n-2^k)$ and $\delta_{a,b}$ is Kronecker delta. In case $k=1$ we use the agreement that $\Pi_{i=a}^{b}W_i$ is the empty word when $b<a$. 
\end{lem}
\pr For a proof see Addendum 1 below. \bx

\begin{cor}\label{const}
The length of any the proper period of $\Eal_n$ on any constituent $K_{m,2^n-1}$ is exactly $2^n-1$. 
\end{cor}

\pr
As it was shown in the previous lemma, the strings $\Eal_n^{2^{n+1}-1}\{K_{m,2^n-1}\}$ are not periodic. Also the $\Eal_n^{2^{n+1}-1}\{K_{m,2^n-1}\}$of length $2^n-1$ of the string $\Eal_n^{2^{n+1}-1}\{K_{m,2^n-1}\}$ coincides with 
\bx

Thus the following holds:
\begin{theo} \label{expon-shifts}
For any collection $\Eal_n$ each element of its spectrum kernel is a divisor of $2^{n+1}-1$.
\end{theo}

This result greatly reduce the upper bound of numbers that lengths of periods for the collections should be multiple of. Indeed, from theorem~\ref{t2} we get the upper bound $2^{2^n-1}-1$ for lengths of proper periods dealing with collection $\Eal_n$. Now it is reduced to $2^{n+1}-1$.

It follows from theorem~\ref{t2}   in cases of odd $n$ collections $\Eal_n$ have periods of all lengths.\footnote{Though lists of proper period lengths is finite. For instance,  
for $n=3$ the position collection is $(1,2,4,8)$. Its proper period lengths are 1,3,15.} 
Here we present lists of lengths of proper periods for an initial segment of even  values of $n$. 
\footnotesize

\begin{center}
 \begin{table}[h]
\caption{All lengths of proper periods for $\Eal_n$ where $n\le 12 $  and is even.}\label{tab2}
\begin{tabular}{|r||l||l|} 
\hline
$n$ &shifts&lengths of proper periods\\
\hline\hline
2 & 1,\,2,\,4 & 7\\
  \hline 
4 & 1,\,2,\,4,\,8,\,16 &31 \\
  \hline   
6 & 1,\,2,\,4,\,8,\,16,\,32,\,64 & 127\\
  \hline   
8& 1,\,2,\,4,\,8,\,16,\,32,\,64,\,128,\,256 & 511,\,73,\,7\\
  \hline   
10 & 1,\,2,\,4,\,8,\,16,\,32,\,64,\,128,\,256,\,512,\,1024 & 2047,\,89,\,23\\
   \hline      
12 & 1,\,2,\,4,\,8,\,16,\,32,\,64,\,128,\,256,\,512,\,1024,\, 2048,\,4096 & 8191\\
\hline
\end{tabular}
\end{table}
\end{center}

\normalsize

Some different collections of shifts might be equivalent in a sense that the strings defined by them have the same periods. 
\begin{theo}\label{c4}
For any positive integer $a$ the collections $C=(x_1,\dots,x_r)$ and $C+a=(x_1+a,\dots,x_r+a)$ have the same periods.
\end{theo}
\pr
This directly follows from two facts. First, the recursion defined by a collection $C$ depends only on elements of $\partial (C)$. And second, $\partial(C)=\partial(C+a)$. 
\bx

Another connection between different collections of positions is described by the next result.

Let $\eta(n,m)$ is the maximal divisor of $n$ s.t. all prime divisors of $\eta(n,m)$ divide $m$. For instance, $\eta(150,20)=50=2\cdot5\cdot5$ because $2,5$ are only primes dividing $20$. 
From the definition it follows that $\eta(n,m)\mid n$ and $\gcd(n,m)=1\implies \eta(n,m)=1$. 

\begin{theo}\label{t4}
Let $C=(x_1,\dots,x_r)$ is a position collection, $B$ is its spectrum kernel, and $a\in\Z^+$. Position collection  $aC=(ax_1,\dots,ax_r)$ has  spectrum kernel $\{b\eta(a,b)|b\in B\}$.
\end{theo}
\pr See Addendum 2 below. \bx

\begin{theo}\label{sym}
Collections $C=(x_1,x_2,\dots,x_{r-1},x_r)$ and  $C^R=(a-x_r,a-x_{r-1},\dots,a-x_1)$ where $a=x_1+x_r$ have the same list of length of periods and thereby the same kernel. 
\end{theo}
\pr
One proof is based on symmetry ``left-right''. Indeed, we can define recursion in the direction from right to left using the equation in a form  
\bea\label{reflect}
\begin{cases}
L(i+x_1-x_r)=L(i+x_1-x_{r-1})+\dots\\
+L(i+x_1-x_2)+L(i+x_1-x_1), i=0,\dots,n-1, \mm2.
\end{cases}
\eea 
This means that if two side infinite string $L$ satisfying the recursion defined by collection $C = (x_1,\dots,x_r)$ satisfies also the recursion~\ref{reflect} defined by $C^R$ in the direction from right to left. Hence the reflection $L^R$ satisfies the recursion defined by the collection $C^R$ but already in the standard direction from left to right. Since all possible periods for strings $L^R$ are reversions of possible periods of strings $L$ and vice versa, the lengths of periods for $C$ and $C^R$ are the same. From here kernels are the same.
\bx 
\begin{ex}{\rm
Collections $(1,5,8,9,11), (1,3,4,7,11)$ have the same list of proper periods lengths: $3,7,21,31,93,217,651$ because $(1,5,8,9,11)^R= (1,3,4,7,11)$.
}\bx
\end{ex}

\ \\

We can introduce a characteristic of reversibility for a collection $C=(x_1,\dots,x_r)$ of shifts as a number  $\mathfrak s=\frac1{1+|\mathcal P|}$ where $\mathcal P$ is the maximal (w.r.t. power) set of relatively prime lengths of periods of the collection. We call this number {\sl reversibility index}. The collections with index 1 contain one shift only, i.e. look like $C=(x_1)$, because $\mathcal P=\emptyset$ for them. All other has stability lesser 1. The collections $C=(x_1,\dots,x_r)$ with even $r$ has zero reversibility index because any prime number is their period. However for odd $r$ collections have finite $|\mathcal P|$. Hence only collections of even powers have zero index. 

\begin{theo}
For any integer $m$ there exist infinitely many collections $C=(x_1,\dots,x_r)$ (even with pairwise different $\partial(C)$ if $m>1$) whose reversibility index is equal to $\frac1{m}$. 
\end{theo}
\pr 
$m$ given consider collection $C=(0,1,2,\dots,a(m)-1)$ where $a(m) = \Pi_{i=1}^{m-1}p_i$ and $\{p_1,\dots,p_{m-1}\}$ is any set of prime numbers s.t.  $2<p_1<p_2<\dots p_{m-1}$. The number $a(m)$ is odd and $a(m)-1$ is even. 
On the other side if $S$ is a finite set of primes and set $\ovl S$ is constructed from all products of primes from $S$, then the power of any system of elements of $\ovl S$ which are relatively prime pairwise does not exceed $|S|$.\footnote{The power of any system of subsets of $S$ that are pairwise disjoint and non empty does not exceed $|S|$. }   Hence $|\Pal|=m-1$ according to theorem~\ref{blocks-t}. Therefore the reversibility index for $C$ is $\frac1m$. Since there exist infinitely many different pairwise sets of $m-1$ primes satisfying the conditions above, the statement is proven.
\bx

Theorem~\ref{t2}  gives an upper limit $2^{x_r-x_1}-1$ for lengths of proper periods a collection $C=(x_1,\dots,x_r)$  might have. As we saw (theorem~\ref{blocks-t}) for block-collection $(x_1,\dots,x_r)$ the upper limit is $x_r-x_1+1$. And for exponential collections it is $2x_r-x_1$, theorem ~\ref{expon-shifts}. We {\sl conjecture} that for collections $(x_1,x_2,x_3)$ where $x_2=x_1+1,x_3=x_1+2^m,m>0,$ the upper limit for proper cycle lengths is $(x_3-x_1)^2-1$.

\newpage

\setcounter{lem}{3}

\section*{Addendum 1: a proof of the lemma 4}\label{exp}
\begin{lem}\label{l7}
The following holds\footnote{For convenience we assume supports of strings to start with number 1.}
\bea\label{form-lem7}
\begin{cases}
\Eal_n^{2^{n+1}-1}\{K_{m,2^n-1}\}=\\
K_{m,2^n-1}\cdot K_{\eta,\eta}\cdot \Pi_{j=k}^{n } 
 \left[\left(\Pi_{i=0}^{k-2}K_{1,2^i}\right)\cdot K_{1,2^j-2^{k-1}+1-m\,\delta_{n,j}}\right]\cdot K_{m,2^{n}-1}
 \end{cases}
\eea
 where ``$\cdot$'' means concatenation of words, $k=\mu z[s_z\ge 2^n-m], \eta=m - (2^n-2^k)$ and $\delta_{a,b}$ is Kronecker delta. In case $k=1$ we use the agreement that $\Pi_{i=a}^{b}W_i$ is the empty word when $b<a$. 
\end{lem}
\pr
For $n\le3$ the lemma  can be checked directly:\\
For $n=2,k=2$ we deal with case $k=n$. When $m=1$ the result is $1001110100$; when $m=2$ we have the resulting word $0100111010$.\\
For $n=3,m=5$ the result is $0^410^21110^2110^510^2$\\
For $n=3,m=6$ the result is $0^510^21110^2110^510$.

Obviously $1\le m< 2^n,\  1\le k\le n$. Then by the definition of $k$ we have $ 2^n-s_{i}\le m\iff  k\le i$. From here 
$m>2^n-2^k$. Therefore $\eta\ge 1$ and $j\ge k\implies 2^j-2^{k-1}+1-m\,\delta_{n,j}\ge1$. The latter follows from $k>i\implies m\le 2^n-2^i$ which in its turn is a consequence of the equivalence $ 2^n-s_{i}\le m\iff  k\le i$ above. So, $m\le 2^n-2^{k-1}$.\footnote{The number of 1(s) depends on $k(m),n$ and is $3+k(n-k+1)$. For $k=1$ and $k=n$ the numbers are equal to $n+3$. For $k=2$ we have $2n+1$ units, etc.}  

The formula from the condition simplifies in cases $k\in\{1,2,n\}$. If $k=1$ then $m=2^n-1$; $k=2$ implies $m\in\{2^n-3, 2^n-2\}$; and when $k=n$ we have $m\le2^{n-1}$. The simplified\footnote{The forms contains single product $\Pi$ on $i$ or $j$ instead of double product on $j,i$ and are directly deducible from the general form.} values of $\Eal_n^{2^{n+1}-1}\{K_{m,2^n-1}\}$  are the following
\bea\label{k=1}
K_{2^n-1,2^n-1}\cdot 1\cdot\left[ \Pi_{j=1}^{n-1 }K_{1,2^j}\right]\cdot 1\cdot K_{2^n-1,2^{n}-1},& k=1, 
\eea
\bea\label{k=n}
K_{m,2^n-1}\cdot K_{m,m}\cdot  
\left(\Pi_{i=0}^{n-2}K_{1,2^i}\right)\cdot K_{1,2^{n-1}+1-m}\cdot K_{m,2^{n}-1},& k=n,
\eea
\bea\label{k=2}
K_{m,2^n-1}\cdot K_{\eta,\eta}\cdot\Pi_{j=2}^{n-1}\left[1\cdot K_{1,2^j-1}\right]\cdot1\cdot  K_{1,2^n-1-m}\cdot K_{m,2^n-1},
& k=2. 
\eea
We consider these cases before the general. Let start with prefixes preceding the product in parentheses. 
In case $k=1$ we have $m=2^n-1$ and the last 1 of the prefix $K_{2^n-1,2^n-1}\cdot 1$ caused by application of $s_1$ to the position $m$ in $K_{2^n-1,2^n-1}$. In case $k=2$ instead of $s_1$ the shift $s_2$ acts while $\eta=m-2^n+4\in\{1,2\}$. In case $k=n$ obviously the position of the last 1 in the prefix is $m+s_k=m+2^n-1$. This explains the suffix $K_{m,m}$ of the prefix. Note that in all the cases $k=1, 2, n$ the considered prefixes have exactly two units with distances $1,3,2^n-1$ respectively between them.

In case $k=1$ the prefix looks like $0^{2^n-2}11$ and therefore the first term $K_{1,2}$ is caused by application of $s_1$ to calculate a symbol on the position $2^n+1$ (which is 1) and shifts $s_1,s_2$ simultaneously for the symbol 0 on the position $2^n+2$. 

Let $\Acal_k$ is a subset of shifts from $\partial\Eal_n$ that act (i.e. contribute 1 into the result) at position $k$. With this term we can say that $\Acal_{2^n+1}=\{s_1\}$ and $\Acal_{2^n+2}=\{s_1,s_3\}$. Thus in case $k=1$ we get a block $1110$ which ends at position $2^n+2$. There is no 1 before the block reachable for shifts from $\Eal_n$ at further positions $x,x>2^n+2$. 

In case $k=n$ two first terms $K_{1,1},K_{1,2}$ of the product (just as above) finish the block $1110$ whose the last position is $2^n+m+2$. This time there is the unit at the position $m$ preceding the block $K_{m,m}\cdot K_{1,1}\cdot K_{1,2}$ but it is not reachable (for shifts from $\Eal_n$) from positions going after $2^n-1+m$. 
Let us rewrite (\ref{k=n}) in form
\bea\label{k=n,a}
K_{m,2^n-1}\cdot K_{m,m}\cdot  \left(\Pi_{i=0}^{n-1}K_{1,2^i}\right)\cdot 
K_{1,2^{n}-m}
\eea
and  prove by induction on $j$ that the common remaining part  $\Pi_{j=2}^{n-1}K_{j,2^i}$ of products from (\ref{k=1}), (\ref{k=n,a}) is the next block of the result of $\Eal_n^{2^{n+1}-1}$ on $K_{2^n-1,2^n-1},K_{m,2^n-1}$. 
For $j=2$ word $W_2$ of length 4 that continues the prefix with $K_{1,4}$ is checkable directly. Let $a$ be a position of the first 1 of the just found sub-word $1110$ in $\Eal_n^{2^{n+1}-1}\{K_{m,2^n-1}\}$. It is easy to check that  $\Acal_{a+4}=\{s_3\}$ implying $W_2=1W'$. Then $\Acal_{a+5}=\{s_1,s_3\}$ and we get $W_2=10W''$. Further, $\Acal_{a+6}=\emptyset,\Acal_{a+7}=\{s_3,s_7\}$. In other words $W_2=1000=K_{1,4}$. 

Assume we have proven the statement for $j=h,h<n-1,$ and consider $j=h+1$.
That is we need to prove that $W_{h+1}=W_h\cdot K_{1,2^{h+1}}$. Clearly the last position of $W_h$ is $a+
\sum_{j=0}^h2^j=a+2^{h+1}-1$. Therefore we should calculate symbols of $W_{h+1}$ at positions from $a+2^{h+1}$ up to $a+2^{h+2}-1$ inclusively. Let start with position $a+2^{h+1}$. In case $k=1$ we can write $\Acal_{a+2^{h+1}}=\{s_i| W_h(a+2^{h+1}-s_i)=1\}$ because no 1 occurs before $W_h$ in $\Eal_n^{2^{n+1}-1}\{K_{n,2^n-1}\}$. In case $k=n$ there is 1 occupying the position $m, m\le 2^{n-1},$ in $\Eal_n^{2^{n+1}-1}\{K_{m,2^n-1}\}$. However as we showed above it is not reachable by available shifts.  Hence the equation $\Acal_{a+2^{h+1}}=\{s_i| W_h(a+2^{h+1}-s_i)=1\}$ fits the case $k=n$ as well. 

By assumption of the induction we have that 1 stays on all positions  $\{a,a+2^t|t=\overline{0,h}\}$ only. Hence we get $(a+2^{h+1}-s_i=a+2^t)\vee (a+2^{h+1}-s_i=a)$. The latter does not work because $s_i$ is odd. The condition reduces to $s_i=2^{h+1}-2^t$ and satisfies at $t=0$ only implying $i={h+1}$. Hence $\Acal_{a+2^{h+1}}=\{s_{h+1}\}$ and $W_{h+1}(a+2^{h+1})=1$.  

We state that 
all other positions from $a+2^{h+1}+1$ up to $a+2^{h+2}-1$ inclusively are filled with zeros.  Let $g$ be an integer s.t. $2^{h+1}+1\le g\le 2^{h+2}-1$. We have $\Acal_{a+g}=\{s_i|W_h(a+g-s_i)=1\}$. Since the position $a+2^{h+1}$ is occupied by one  we have 
\beas
s_i\in \Acal_{a+g}\iff\\
\left(\exists t,0\le t\le h+1\right) \left[a+g-s_i=a+2^t\right]\vee (a+g-s_i=a) \iff\\
\left(\exists t,0\le t\le h+1\right)\left[g=s_i+2^t\right]\vee (g=s_i)
\eeas
From the first term of the disjunction (the first disjunct)\footnote{We will use the word ``disjunct'' for terms of this disjunction and disjunction obtained below in this proof.}: since $0\le t$ and $g,t<2^{h+2}$ it holds $i\le h+1$. On the other hand, $g>2^{h+1}$. Hence $\max\{i,t\}\ge h+1, \min\{i,t\}< h+1$. The second disjunct is fulfilled when $i=h+2$ only due to the restrictions $g$ satisfies. Hence $\Acal_{a+g}\neq\emptyset \iff g\in\{2^{h+1}+2^t-1| t=0,\dots,h+1\}$. Now, if $g=2^{h+1}+2^t-1, t< h+1,$ obviously $\Acal_{a+g}=\{s_t,s_{h+1}\}$. If $g=2^{h+2}-1$ then $\Acal_{g+a}=\{s_{h+1},s_{h+2}\}$ because $s_{h+2}\in\Eal_n$ yet: $h+2\le n$. The induction on $j$ is over.

Now we need to prove that the remaining suffixes of the results of action of $\Eal_n^{2^{n+1}-1}$ on $K_{2^n-1,2^n-1}, K_{m,2^n-1}$ are respectively $1\cdot K_{2^n-1,2^n-1}$ and  $K_{1,2^n-m}$ (see (\ref{k=n,a})).

Since $1\cdot K_{2^n-1,2^n-1}=K_{1,2^n-1}\cdot1$ we can apply the consideration from the induction step above (getting case $h=n-1$). The difference is only that for the last position $b$ of the suffix $K_{1,2^n-1}\cdot1$ we have $\Acal_b=\{s_n\}$ because $s_{n+1}\notin\Eal_n$. This is why the position $b$ is occupied by 1, not 0 as it took place on the induction step above. 

Let $b$ be the first position in the conjectured suffix $K_{1,2^n-m}$. It is easy to calculate that $b=2^n-1+m+\sum_{i=0}^{n-1}2^i+1=m+2^{n+1}-1$. The positions of units preceding the unit on the position $b$ in the word (\ref{k=n}) are: $m,m+2^n-1,m+2^n-1+2^0,\dots,m+2^n-1+2^{n-1}$. So, positions $m,m+2^n-1$ are not reachable from position $b$ by any $s_i\in\Eal_n$. Other positions before $b$ can be written as $m+2^n-1+2^t,t=\overline{0,n-1}$.
The equation $m+2^{n+1}-1-(m+2^n-1+2^t)=2^{n}-2^t\in \Eal_n$ has only solution $t=0$. Hence $\Acal_b=\{s_n\}$ and position $b$ is occupied by 1 indeed. 

Now let us find $\Acal_{b+h},0<h<2^n-m$ assuming that for any $h'$ such that $0<h'<h$ the position $b+h'$ is occupied by 0. The list of units on positions before $b+h$ to take into account is\footnote{Positions $m,m+2^n-1$ are not reachable even from $b$, see above.}   
\beas
m+2^n-1+2^0,\dots,m+2^n-1+2^{n-1},m+2^{n+1}-1,
\eeas
and the equation for $\Acal_{b'}$ is 
\beas
(h\in\Eal_n)\vee\exists t\in\{0,\dots,n-1\}\left[m+2^{n+1}-1+h-(m+2^n-1+2^t) \in \Eal_n\right].
\eeas
The first disjunct with $m\ge 1,h<2^n-m$ can be rewritten as $h\in\{s_1,\dots,s_{n-1}\}$. The condition $2^n+h-2^t\in\Eal_n$  from the second disjunct can be rewritten as $\exists i[h=2^i+2^t-2^n-1]$. Because $t\le n-1$ from $h>0$ it follows that the unique solution is $i=n$ and $h=2^t-1$. Due to the latter equality is in accord with the first disjunct we have that $\Acal_{b+h}=\{s_t,s_n\}$ if $h=s_t$ and $\Acal_{b+h}=\emptyset$ otherwise. As $t<n$ we have that $|\Acal_{b+h}|$ is even and all positions $b+h$ are occupied by zeros. 

Thus the formulas (\ref{k=1}) and (\ref{k=n}) are correct. 

In case $k=2$ we have to prove (\ref{k=2}):
\beas
K_{m,2^n-1}\cdot K_{\eta,\eta}\cdot\Pi_{j=2}^{n-1}\left[1\cdot K_{1,2^j-1}\right]\cdot1\cdot  K_{1,2^n-1-m}\cdot K_{m,2^n-1}.
\eeas
It is easy to check that $(2^n-1+\eta)-m=s_k=s_2=3$. Also directly one can check that the first factor $1\cdot K_{1,3}$ concatenated to $K_{m,2^n-1}\cdot K_{\eta,\eta}$ is a prefix of $\Eal_n ^{2^{n+1}-1}{K_{m,2^n-1}}$. The list of positions of units in $K_{m,2^n-1}\cdot K_{\eta,\eta}\cdot 1\cdot K_{1,3}$ is $m,m+3,m+4,m+5$. In particular, $\Acal_{m+6}=\{s_1,s_2\},\Acal_{m+7}=\{s_2,s_3\}$. We prove that product $\Pi_{j=2}^{n-1}\left[1\cdot K_{1,2^j-1}\right]$ correctly describes the continuation by induction on $j,2\le j<n$. From the induction hypothesis the list of units for prefix $K_{m,2^n-1}\cdot K_{\eta,\eta}\cdot\Pi_{j=2}^{h}\left[1\cdot K_{1,2^j-1}\right], 2\le h<n-1,$ is
\bea\label{ones list}
m,m+3,m+4,m+5,\dots, m+2^h,m+2^h+1
\eea
whereas the last symbol occupies the position $m+2^h+2^h-1=m+2^{h+1}-1$. We can define $\Acal_{m+2^{h+1}}$ 
by the condition 
\beas
s_i\in \Acal_{m+2^{h+1}}\iff \exists t, 2\le t\le h[s_i=2^{h+1}-2^t \vee s_i=2^{h+1}-1-2^t]
\eeas
since the condition $(s_i=2^{h+1})\vee(s_i=2^{h+1}-3)$ can be omitted because $n\ge 4\implies h+1\ge3$. It is clear that the condition $s_i=2^{h+1}-2^t$ is impossible when $t\ge 2$. The remaining condition has a unique solution $t=i=h$. Therefore position $2^{h+1}$ is occupied by 1. 

Obviously $s_1\in \Acal_{m+2^{h+1}+1}$. If $i>1$ then as above for $\Acal_{m+2^{h+1}+1}$ we arrive at condition
\beas
s_i\in \Acal_{m+2^{h+1}}\iff \exists t, 2\le t\le h[s_i=2^{h+1}+1-2^t \vee s_i=2^{h+1}-2^t], i>2.
\eeas
Both are not solvable because $t>1$ and $2^{h+1}-2^t$ is even. Hence position $m+2^{h+1}+1$ is occupied by 1 again. 

Let us prove that the next $2^{h+1}-2$ positions are occupied by zeros. Let assume $1\le q\le2^{h+1}-2$ and prove that if for any $q',0<q'<q$ position $m+2^{h+1}+1+q'$ is occupied by zero then the same is true for position $m+2^{h+1}+1+q$.
We have 
\bea\label{cond-k=2}
\hspace{3em}\begin{cases}
s_i\in\Acal_{m+2^{h+1}+1+q}\iff  \big\{(\exists t\in\{2,\dots, h+1\})\big[s_i=2^{h+1}+1+q-2^t\quad\vee\\  \qquad  \quad s_i=2^{h+1}+q-2^t\big]\big\}\vee \left(s_i=2^{h+1}+1+q\right) \vee \left( s_i=2^{h+1}+1+q-3\right).
\end{cases}
\eea
There are 4 disjuncts in the condition:\\
(i) $\exists t, 2\le t\le h+1\left[s_i=2^{h+1}+1+q-2^t\right]$;\\
(ii) $\exists t, 2\le t\le h+1\left[s_i=2^{h+1}+q-2^t\right]$;\\
(iii) $s_i=2^{h+1}+q+1$;\\
(iv) $s_i=2^{h+1}+q-2$.\\
Clearly no two of these conditions are compatible. Let us analyse them one by one. 

For (i) on one hand we have $s_i<2^{h+2}-4$ due to $t\ge 2$ and given bounds to $q$. This entails $i\le h+1$. On the other hand, from $q=2^i+2^t-(2^{h+1}+2)>0$ we have $2^i+2^t\ge 2^{h+1}+3$. Hence if $t<h+1$ then $i=h+1$ and hence $q=2^t-2$. Thus that is one series of solutions to (i): $i=h+1,q=2^t-2, 1<t<h+1$. Another series arises when we set $t=h+1$. In this case we have $q=2^i-2,i=2,\dots,h+1, t=h+1$.

Similarly\footnote{On one hand we have $s_i<2^{h+2}-5$ due to $t\ge 2$ and given bounds to $q$. This entails $i\le h+1$. On the other hand, from $q=2^i+2^t-(2^{h+1}+1)>0$ we have $2^i+2^t\ge 2^{h+1}+2$. Hence if $t<h+1$ then $i=h+1$ and hence $q=2^t-1$. Thus that is one series of solutions to (i): $i=h+1,q=2^t-1, 1<t<h+1$. Another series arises when we set $t=h+1$. In this case we have $q=s_i,i=1,\dots,h, t=h+1$.} for (ii) we get one series $ q=2^t-1,1<t<h+1,i=h+1,$  and another one $q=s_i, i= 1,\dots,h,t=h+1$.

For (iii) the unique solution is $q=2^{h+1}-2, i=h+2$. Indeed, if $q=s_i-(2^{h+1}+1)$ then because $q>0$ it holds $i\ge h+2$. However $i> h+2$ is incompatible with the upper bound $2^{h+1}-2$ for $q$.

Similarly from (iv) $i=h+1$ follows and therefore we get unique solution  $q=1$. 

From here $\Acal_{m+2^{h+1}+1+q}\neq\emptyset\iff \underset{f=1,\dots,h+1}{\vee}(q=2^f-1 \vee(q=2^f-2)$. Combining solutions to (ii) and (iv) we get $q=2^f-1\implies \Acal_{m+2^{h+1}+1+q}=\{s_f,s_{h+1}\}$; combining solutions to (i), (iii) we arrive at $q=2^f-2\implies \Acal_{m+2^{h+1}+1+q}=\{s_f,s_{h+1}\}, t=2,\dots,h,$ and if $f=h+1$ then $\Acal_{m+2^{h+1}+1+q}=\{s_f,s_{h+2}\}$.
In other words, all $|\Acal_{m+2^{h+1}+1+q}|$ are even if $q=\overline{1,2^{h+1}-2}$. This finishes the induction on $j$.

It remains only to prove the correctness of the suffix $1\cdot K_{1,2^n-1-m}\cdot K_{m,2^n-1}$. We rewrite it as $1\cdot K_{1,2^n-2}\cdot K_{1,2^n-m}$ and first prolong the consideration of correctness the product $\Pi_{j=2}^{n-1}\left[1\cdot K_{1,2^j-1}\right]$ from (\ref{k=2}) given above. The block $1\cdot K_{1,2^n-2}$ corresponds to step $j=n$ and we need only to explain the difference between the length $2^n$ of the generic term $1\cdot K_{1,2^j-1}$ for $j=n$ and the length $2^n-1$ of the block $1\cdot K_{1,2^n-2}$. The conditions (\ref{cond-k=2}) are applicable here as well but now $h+1=n$ and for the position $\bar q=2^n-2$ the equality $\Acal_{m+2^{h+1}+1+\bar q}=\{s_{h+1},s_{h+2}\}$ looks like $\Acal_{m+2^{n}+1+\bar q}=\{s_{n},s_{n+1}\}$. And because $s_{n+1}\notin\Eal_n$ we get that in fact $\Acal_{m+2^{n}+1+\bar q}=\{s_{n}\}$. The latter enforces 1 on the position $m+2^{n}+1+\bar q=m+2^n+1+2^n-2=m+2^{n+1}-1$. This position is exactly the first position of the last block $K_{1,2^n-m}$.\footnote{
Indeed, according to the formula (\ref{k=2}) the position of the last 1 is $2^n-1+\eta+2^n-4+1+2^n-1-m+m= 2^n-1+m-2^n+4+2^n-4 +1+2^n-1-m+m= m+2^{n+1}-1$.}

To finish the case $k=2$ we need only to explain why the suffix of the (\ref{k=2}) of length $2^n-m-1$ consists of zeros. 
We know now positions of all units preceding the suffix. So we can complete the list (\ref{ones list}) as follows:
\beas
m,m+3,m+4,m+5,\dots, m+2^h,m+2^h+1,\dots,m+2^n,m+2^n+1,m+2^{n+1}-1.
\eeas
Since from a position $m+x,x\ge 2^{n+1},$ by shifts from $\Eal_n$ are reachable only maybe units on the last three positions $m+2^n,m+2^n+1,m+2^{n+1}-1$ the equation for $\Acal_{x}$ is
\beas
s_i\in\Acal_x\iff (x-s_i=2^n) \vee (x-s_i=2^n+1) \vee (x-s_i=2^{n+1}-1).
\eeas
On the other hand, because of $k=2$ , the upper bound for $x$ is $2^{n+1}-1+(2^n-1)-m-1$ which does not exceed $2^{n+1}-1+s_2-1$. In other words $x\le2^{n+1}+1$.\footnote{The number of zeros after the last 1 does depend on one of two possible positions of $m$ for case $k=2$ but not exceed 2.} From here the first disjunct in the equation does not work and other two have the only solution for $s_i\in\Acal_x=\{s_1,s_n\}$ if $x=1$ and $\empty$ otherwise. 

\ \\

Let us try a general proof. First let us rewrite the formula from the condition in the form
\bea\label{new-view}
\begin{cases}
\Eal_n^{2^{n+1}-1}\{K_{m,2^n-1}\}=\\
K_{m,2^n-1}\cdot K_{\eta,\eta}\cdot \Pi_{j=k}^{n-1 } 
 \left[W_k\cdot K_{1,2^j-2^{k-1}+1}\right]\cdot W_k\cdot K_{1,2^n-2^{k-1}}\cdot K_{1,2^{n}-m}
\end{cases}
\eea
where $W_k=\Pi_{i=0}^{k-2}K_{1,2^i}$. Now, $2^n-1+\eta$ is the first position from which the unit on position m is reachable due to the definition of $k$ and equality $2^n-1+\eta=s_k+m$. This explains the sub-word $K_{\eta,\eta}$ succeeding $K_{m,2^n-1}$. Thus we deal with a word $K_{m,2^n-1}\cdot K_{\eta,\eta}$ having two units only on positions $m$ and $m+s_k$. 

Next step is proving the stated structure of the block 
\bea\label{products}
\Pi_{j=k}^{n-1 } \left[W_k\cdot K_{1,2^j-2^{k-1}+1}\right]\cdot W_k\cdot K_{1,2^n-2^{k-1}}
\eea
by an induction on $j$ with inner induction on $i$. The part $W_k\cdot K_{1,2^n-2^{k-1}}$ corresponds to the value of  generic factor $W_k\cdot K_{1,2^j-2^{k-1}+1}$ at $j=n$ with only difference in the length of  $K_{1,2^n-2^{k-1}}$ which must be explained. 

Since cases $k=n,k\le2$ were considered above we assume that $2<k<n$.  Therefore the base of the external induction consists of a proof of the following structure of the first factor in the external product
\bea\label{bas-two-prods}
\Pi_{i=0}^{k-2}K_{1,2^i}\cdot K_{1, 2^k-2^{k-1}+ 1}.
\eea

The base of the inner induction for (\ref{bas-two-prods}) states that $K_{1,2^0}$ is a prefix of the word succeeding $K_{m,2^n-1}\cdot K_{\eta,\eta}$. Since $K_{1,2^0}=1$ this is obvious - this unit is generated by shift $s_1$ on the last unit of $K_{m,2^n-1}\cdot K_{\eta,\eta}$.

Now we assume $K_{m,2^n-1}\cdot K_{\eta,\eta}\cdot \Pi_{i=0}^{h}K_{1,2^i}, h<k-2,$ is a prefix of 
\bea\label{RHS}
K_{m,2^n-1}\cdot K_{\eta,\eta}\cdot \Pi_{j=k}^{n-1 }  \left[W_k\cdot K_{1,2^j-2^{k-1}+1}\right]\cdot W_k\cdot K_{1,2^n-2^{k-1}}\cdot K_{1,2^{n}-m}
\eea
and prove the same for $K_{m,2^n-1}\cdot K_{\eta,\eta}\cdot \Pi_{i=0}^{h+1}K_{1,2^i}$. For that we write out positions of units in $K_{m,2^n-1}\cdot K_{\eta,\eta}\cdot \Pi_{i=0}^{h}K_{1,2^i}$: 
\bea\label{list-bas}
m,m+s_k,m+s_k+2^0,\dots, m+s_k+2^h, 
\eea
wereas all other positions to $m+s_k+2^{h}+2^{h}-1$ inclusively are occupied with zeros. 
Since $h+1\le k-2$ the next $K$-block in (\ref{RHS}) should be $K_{1,2^{h+1}}$. The distances from position $m+s_k+2^{h+1}$ to the units in the list (\ref{list-bas}) are
\beas
s_k+2^{h+1},2^{h+1},2^{h+1}-2^0,\dots, 2^{h+1}-2^h.
\eeas
Since $h>0$ all numbers in the list are even apart from $s_k+2^{h+1},2^{h+1}-2^0$. Equation $2^i-1=2^k-1+2^{h+1}$ or $2^i=2^k+2^{h+1}$ has a solution $i=h+2$ only if $k=h+1$. But $h<k-2$. Therefore the only term of the above list that belongs to $\Eal_n$ is $2^{h+1}-2^0$. Hence $|\Acal _{m+s_k+2^{h+1}}|=1$ and this position is occupied by 1. 
So we prolong the list (\ref{list-bas}) as
\bea\label{1-pos} 
m,m+s_k,m+s_k+2^0,\dots, m+s_k+2^h, m+s_k+2^{h+1}
\eea 
and prove that all  positions $q+m+s_k,q=2^{h+1}+1,\dots,2^{h+2}-1$ are occupied by zeros. It holds 
\bea\label{condsWk}
\hspace{2em} s_i\in\Acal_{m+s_k+q}\iff\exists t\in\{0,\dots,h+1\}[s_i=q-2^t]\vee(s_i=q)\vee(s_i=s_k+q).
\eea
Any equation $q=s_i+2^t$ has two solutions when $q$ is in the considered segment. One is defined by $i=h+1,0\le t\le h$. Another is defined by $t=h+1,i=1,\dots,h+1$. 

There is one only number of kind $s_i$ in segment $[2^{h+1}+1,2^{h+2}-1]$. It is $2^{h+2}-1$. 

Equation $2^i=2^k+q$ implies $q=2^k, i=k+1$. But $q<2^{h+2}\le 2^k$. Therefore the last disjunct in (\ref{condsWk}) does not work for the positions we consider now. 

As a result we have $\Acal_{m+s_k+q}\neq\emptyset\iff q=2^{h+1}+s_t,t=0,\dots,h+1$. Also  
\beas
\Acal_{m+s_k+q}=
\begin{cases}
\{s_t,s_{h+1}\},& q=2^{h+1}+s_t,t<h+1,\\ 
\{s_{h+1},s_{h+2}\},& q=2^{h+2}-1.
\end{cases}
\eeas
The induction on $i$ is closed. Now we prove that the $K$-block $K_{1, 2^k-2^{k-1}+ 1}$ from (\ref{bas-two-prods}) is correct continuation of the product before it.  In other words we need to prove that position\footnote{$h+1=k-2$}
$m+s_k+2^{k-1}$ is occupied by 1 and next $2^{k-1}$ positions are occupied by zeros. The distances from position  $m+s_k+2^{k-1}$ to positions from the list (\ref{1-pos}) where $h=k-3$ are 
\beas
s_k+2^{k-1},2^{k-1},2^{k-1}-1,2^{k-1}-2^{1},\dots, 2^{k-1}-2^{k-3},2^{k-1}-2^{k-2}.
\eeas
Among two odd distances $s_k+2^{k-1},2^{k-1}-1$ the first is not element of $\Eal_n$. The second however is. Therefore the position $m+s_k+2^{k-1}$ is occupied by 1 indeed. For the next $2^{k-1}$ positions $q+m+s_k,2^{k-1}<q\le 2^k,$ we need slightly change the conditions (\ref{condsWk}) as follows
\bea\label{condsafterW}
\hspace{2em} s_i\in\Acal_{m+s_k+q}\iff\exists t\in\{0,\dots,k-1\}[s_i=q-2^t]\vee(s_i=q)\vee(s_i=s_k+q) 
\eea
because now we base on the list 
\beas
m,m+s_k,m+s_k+2^0,\dots, m+s_k+2^t,\dots, m+s_k+2^{k-1}
\eeas
of units before positions $q+m+s_k$. 

Any equation $q=s_i+2^t$ has two solutions when $q$ is in the segment $[2^{k-1}+1,2^k-2]$. One is defined by $i=k-1,1\le t < k-1$. Another is defined by $t=k-1,i=1,\dots,k-2$. 
There is one solution $i=t=k-1$ for $q=2^{k}-1$ and one solution  $t=0,i=k,q=2^k$ for $q=2^k$. 
There is one only number of kind $s_i$ in segment $[2^{k-1}+1,2^{k}]$. It is $2^{k}-1$. 
Equation $2^i=2^k+q$ implies $q=2^k, i=k+1$.   

As a result we have $\Acal_{m+s_k+q}\neq\emptyset\iff q=2^{h+1}+s_t,t=0,\dots,h+1$. Also  
\beas
\Acal_{m+s_k+q}=
\begin{cases}
\{s_t,s_{k-1}\},& q=2^{k-1}+2^t-1,1\le t< k-1,\\ 
\{s_{k-1},s_{k}\},& q=2^{k}-1,\\
\{s_k,s_{k+1}\}, & q=2^k.
\end{cases}
\eeas
Thus the basis of the induction on $j$ is completed. 

\ \\

The scheme of transition $j\to j+1$ is the same as above: we write out the list of  positions of units, and compile the conditions of reachability by shifts from $\Eal_n$ any of the units from a position we currently consider. Solving the conditions we install what symbol from $\{0,1\}$ occupies the current position. 

So, assume the correctness of the prefix $K_{m,2^n-1}\cdot K_{\eta,\eta}\cdot \Pi_{j=k}^{h } 
\left[W_k\cdot K_{1,2^j-2^{k-1}+1}\right],h<n-1,$ in (\ref{new-view}) is proved. From the prefix inspection, the list of positions of units is 
\beas
\begin{cases}
m,m+s_k,m+s_k+2^0,\dots, m+s_k+2^{k-1},\dots,\\
m+s_k+\sum_{j=k}^t2^j+2^0,\dots m+s_k+\sum_{j=k}^t2^j+2^{k-1},\dots\\
m+s_k+\sum_{j=k}^h2^j+2^0,\dots m+s_k+\sum_{j=k}^h2^j+2^{k-1}
\end{cases}
\eeas
or
\beas
m,m+s_k, 
a_0+2^0,\dots, a_0+2^{k-1},\dots,
a_t+2^0,\dots, a_t+2^{k-1},\dots, 
a_h+2^0,\dots,a_h+2^{k-1}
\eeas
where $a_0=m+s_k$ and $a_t=m+s_k+\sum^{t}_{j=k} 2^j=m+s_k+2^{t+1}-2^{k}=m+s_{t+1}$ where $t$ runs over $k,\dots,h$. So we can rewrite the list as follows 
\bea\label{gen-case-list}
\begin{cases}
m,m+s_k,\\
m+s_k+2^0,\dots, m+s_k+2^u,\dots,  m+s_k+2^{k-1},\\
\dots\\
m+s_{t}+2^0,\dots,m+s_t+2^u,\dots,  m+s_{t}+2^{k-1},\\
\dots\\
m+s_{h}+2^0,\dots,m+s_h+2^u,\dots,  m+s_{h}+2^{k-1},
\end{cases}
\eea
where $k\le t\le h$. Also it is easy to see that the last position of the prefix $P_h=K_{m,2^n-1}\cdot K_{\eta,\eta}\cdot \Pi_{j=k}^{h } \left[W_k\cdot  K_{1,2^j-2^{k -1}+1}\right]$. It is $m+s_{h+1}$.\footnote{It is a sum of position $m+s_h+2^{k-1}$ of the last 1  in $P_h$ and the number $2^h-2^{k-1}$ of zeros in $K_{1,2^j-2^{k -1}+1}$.} 

Now we check if the block $W_k\cdot  K_{1,2^{h+1}-2^{k-1}+1}$ is the correct continuation of $P_h$. 
The condition determining $\Acal_{m+s_{h+1}+2^0}$ is:
\beas
s_i\in \Acal_{m+s_{h+1}+2^0}\iff (s_i=s_{h+1}+2^0)\vee(s_i=s_{h+1}+2^0-s_k)\vee\\
\exists t\in\{k,\dots,h\}\exists u\in\{0,\dots,k-1\}[s_i=s_{h+1}+2^0-s_t-2^u].
\eeas

The first disjunct can be cancelled because $s_{h+1}+1$ is an even integer but $s_i$ is odd. 

The second disjunct has no solution when $h\ge k>2$ because $2^{h+1}-2^k\ge 2^k>2$.\footnote{If $2^i=(2^{h+1}-2^k)+2$ then to complete $(2^{h+1}-2^k)$ to a power of 2 the amount $(2^{h+1}-2^k)$  must be equal to 2.}

For 3rd disjunct we can rewrite equation as $2^i=2^{h+1}+2 - (2^t+2^u)$. If $u=0$ then $s_{h+1}-s_t$ must be odd but that is impossible due to $h+1>t\ge k>2$. Let $u>1$. Setting maximal possible values $t:=h,k:=h,u=h-1$ we get $2^i=2^{h+1}+2-2^h-2^{h-1}\ge 2^{h-1}+2$. If we set $t:=k,u:=1$ we get $2^i\le 2^{h+1}+2-(2^k+2^1)=2^{h+1}-2^k\le 2^{h+1}-8$. In the segment $[2^{h-1}+2,2^{h+1}-8]$ one only power of 2 exists: $2^h$. Hence $i=h$ and therefore 
$2^{h+1}-2^i=2^h=2^t+2^u-2$ where $u<k\le t\le h$. From here $t=h,u=1$. No other solution exists. This implies $|\Acal_{m+s_{h+1}+2^0}|=1$ and the position is occupied by 1.

Being still within block $W_k$ from block $W_k\cdot  K_{1,2^{h+1}-2^{k-1}+1}$, we assume that a prefix $W_k\big|_{m+s_{h+1}+1}^{m+s_{h+1}+2^{u+1}-1}$ of $W_k$ is proven to be correct. Now we do next step and prolong the proof on positions $m+s_{h+1}+2^{u+1},m+s_{h+1}+2^{u+1}+q,q=1,\dots, 2^{u+1}-1$. Notice that to be within $W_k$ the number $u+1$ must be lesser $k-1$.
The list of positions of preceding units is\footnote{In each line of the list below the last 1 belongs to a $K$-block following the previous $W_k$-block; lines correspond to parameter $j$ running $k,\dots,t,\dots,h,h+1$.}
\bea\label{list-u}
\begin{cases}
m,m+s_k,\\
m+s_k+2^0,\dots, m+s_k+2^u,\dots,  m+s_k+2^{k-1},\\
\dots\\
m+s_{t}+2^0,\dots,m+s_t+2^u,\dots,  m+s_{t}+2^{k-1},\\
\dots\\
m+s_{h}+2^0,\dots,m+s_h+2^u,\dots,  m+s_{h}+2^{k-1},\\
m+s_{h+1}+2^0,\dots,m+s_{h+1}+2^u
\end{cases}
\eea
and therefore the equation for $\Acal_{m+s_{h+1}+2^{u+1}}$ differs from the previous: 
\bea\label{pos2}
s_i\in \Acal_{m+s_{h+1}+2^{u+1}}\iff (s_i=s_{h+1}+2^{u+1}) \vee(s_i=s_{h+1}+2^{u+1}-s_k)\vee\notag\\
\exists t\in\{k,\dots,h\}\exists v\in\{0,\dots,k-1\}[s_i=s_{h+1}+2^{u+1}-s_t-2^v]\vee\\
\exists v\in\{0,\dots,u\}[s_i=2^{u+1}-2^v].\notag
\eea
Here $0<u+1<k-1< t\le h$ 
and therefore the first disjunct does not have solutions. The second disjunct has no solution either because $s_i$ is odd but $s_{h+1}+2^{u+1}-s_k$ is even.\footnote{$k>2$.} To avoid the problem with parity and fulfil $s_i=s_{h+1}+2^{u+1}-s_t-2^v$ one need to set $v=0$. So, we can write it as $s_i=(s_{h+1}-s_t)+s_{u+1}$. The only opportunity to satisfy the equation is $t=h=u+1$ (in which case $i=u+2=h+1$). However that is impossible because $u+1<h$. The forth disjunct however has the unique solution $s_i=s_{u+1}$. That is why $ \Acal_{m+s_{h+1}+2^{u+1}}=\{s_{u+1}\}$ and the position is occupied by 1. 

Now assuming that for all $q',1\le q'<q$ where $1\le q\le2^{u+1}-1$ a position $m+s_{h+1}+2^{u+1}+q'$ is occupied by 0, we prove that the position $m+s_{h+1}+2^{u+1}+q$ is occupied by 0 as well. 
The list of unit positions is pretty much similar to (\ref{list-u}). The only difference is in the last line where $u$ is replaced with $u+1$, and therefore the equation for $\Acal_{m+s_{h+1}+2^{u+1}+q}$ is:
\bea\label{condsW}
s_i\in \Acal_{m+s_{h+1}+2^{u+1}+q}\iff \hspace{7.5cm}\\
\begin{cases}\notag
 (s_i=s_{h+1}+2^{u+1}+q) \vee(s_i=s_{h+1}+2^{u+1}+q-s_k)\vee\\
\exists t\in\{k,\dots,h\}\exists v\in\{0,\dots,k-1\}[s_i=s_{h+1}+2^{u+1}+q-s_t-2^v]\vee\\
\exists v\in\{0,\dots,u+1\}[s_i=2^{u+1}+q-2^v].
\end{cases}
\eea
Since $u+1<k-1<h<h+1$ 
and $1\le q<2^{u+1}$ we have $2^{h+1}+2^{u+1}<2^{h+1}+2^{u+1}+q< 2^{h+1}+2^h$. No power of 2 exists in the segment $[2^{h+1}+2^{h-1}, 2^{h+1}+2^h],h>2,$ and therefore no solution to the first disjunct exists. 

The second disjunct has no solutions. Indeed, due to $1\le u+1\le k-2, \ 1\le q\le 2^{u+1}-1,\ k\le h$ we have  $2< 2^{u+1}+q\le s_{k-1}<2^h$. Hence  $2^{h+1}-2^{k-1}<2^{h+1}-(s_k-2^{u+1}-q)< 2^{h+1}$ resulting impossibility for $2^{h+1}-(s_k-2^{u+1}-q)$ to be a power $2^i$.

The third disjunct has solutions $t=h,v=u+2,q=2^{u+1}-1,i=h$ taking into account $u+2\le k-1$. No other solution exists.\footnote{The case $q=s_v,t=u+1,s_i=s_{h+1}$ is not possible because $t\ge k>u+1$.}. To see that let first prove that $i=h$ has no alternative. Indeed, we can rewrite the matrix of the 3rd disjunct as $2^i=2^{h+1}+2^{u+1}+(q+1)-2^t-2^v$
and estimate maximal and minimal values of the sum $g=2^{u+1}+(q+1)-2^t-2^v$ having $1\le u+1<k-1,\  2\le q+1\le 2^{u+1},\ k\le t\le h,\  0\le v\le k-1,\ 3\le k\le h$. The maximal value of $g$ is $-2^{k-1}-1$ whereas the minimal value is $-2^h-2^{k-1}+4$. That yields $2^{h-1}+4\le 2^{h}-2^{k-1}+4\le 2^{h+1}+g\le 2^{h+1}-2^{k-1}-1$. So one only power of 2 exists in this interval and it is $2^h$ Hence $i=h$. To get that $t$ must be equal to $h$ and we arrive at the equality $2^v=2^{u+1}+q+1$. From here $q=s_{u+1},v=u+2$. And because of $v\le k-1$ the said is true when $u+1\le k-2$.

The forth disjunct has two series of solutions. One is $q=s_v,i=u+1, v=1,\dots, u+1$. Another one is $v=u+1,s_i=q$ where $q$ runs over list $1,\dots,s_u$. The case $v=u+1,s_i=q=s_{u+1}$ coincides with the last solution in the previous series.\footnote{The case $v=0$ is impossible because $q<2^{u+1}$.} Since $s_i\le s_{u+1}$ and all such $s_i$ are included into the solution, the consideration is complete.  

The summary is that $\Acal_{m+s_{h+1}+2^{u+1}+q}=\emptyset$ if $q\notin\Eal_n$ and defined as follows
\beas
\Acal_{m+s_{h+1}+2^{u+1}+q}=
\begin{cases}
\{s_i,s_{u+1}\},& q=s_1,\dots, s_u,\\
\{s_{u+1},s_h\},& q=s_{u+1}, 
\end{cases}
\eeas
otherwise. Thus all positions $m+s_{h+1}+2^{u+1}+q,q=s_1,\dots,s_{u+1}$ are occupied by zeros. This conclusion finishes 
the induction on $u$ and $W_k$ is passed. \\

Let us consider the first position after the considered block $W_k$. It is $m+s_{h+1}+2^{k-1}$. The equation for $\Acal_{m+s_{h+1}+2^{k-1}}$ is similar to (\ref{pos2})
\beas
s_i\in \Acal_{m+s_{h+1}+2^{k-1}}\iff (s_i=s_{h+1}+2^{k-1}) \vee(s_i=s_{h+1}+2^{k-1}-s_k)\vee\\
\exists t\in\{k,\dots,h\}\exists v\in\{0,\dots,k-1\}[s_i=s_{h+1}+2^{k-1}-s_t-2^v]\vee\\
\exists v\in\{0,\dots,k-2\}[s_i=2^{k-1}-2^v].
\eeas
Due to $2^{h+1}+ 2^{k-1}$ is not a power of 2 the first disjunct fails. The RHS of the second is even whereas the LHS must be odd. To make the RHS in the third disjunct odd we must set $v=0$. So, $i$ should satisfy the inequalities $2^{h}+2^{k-1}\le 2^i\le 2^{h+1}-2^{k-1}$. Clearly that is impossible. Finally in the 4th disjunct we must first set $v=0$. hence the unique solution is $s_i=s_{k-1}$. 

Thus the position is occupied by 1 and according to the structure $K_{1,2^{h+1}-2^{k-1}+1}$ of the sub-word we are considering now, we need to prove that zeros stay on next $2^{h+1}-2^{k-1}$ positions. Let $1\le q\le 2^{h+1}-2^{k-1}$. The conditions for $\Acal_{m+s_{h+1}+2^{k-1}+q}$ are similar to (\ref{condsW}): 
\beas
s_i\in \Acal_{m+s_{h+1}+2^{k-1}+q}\iff \hspace{7.5cm}\\
\begin{cases}\notag
 (s_i=s_{h+1}+2^{k-1}+q) \vee(s_i=s_{h+1}+2^{k-1}+q-s_k)\vee\\
\exists t\in\{k,\dots,h\}\exists v\in\{0,\dots,k-1\}[s_i=s_{h+1}+2^{k-1}+q-s_t-2^v]\vee\\
\exists v\in\{0,\dots,k-1\}[s_i=2^{k-1}+q-2^v].
\end{cases}
\eeas
The only way to satisfy $s_i=s_{h+1}+x, x>0,$ is to put $x=2^{h+1}$. Then $s_i=s_{h+2}$ and $q=2^{h+1}-2^{k-1}$.

The second disjunct can be rewritten as $s_i=s_{h+1}+q-s_{k-1}$. Since $ 2-2^{k-1}\le  q-s_{k-1}\le 2^{h+1}-s_k$ the only way to get a number of sort $s_i$ for some $i$ from $s_{h+1}+q-s_{k-1}$ is to set $q=s_{k-1}$ resulting $s_i=s_{h+1}$. 

Represent the 3rd disjunct in form $s_i=(s_{h+1}-s_t)+(2^{k-1}-2^v)+q$. The restrictions $0\le v\le k-1,\ k\le t  \le h$ imply $s_h\le (s_{h+1}-s_t) \le s_{h+1}-s_k,\ 0\le (2^{k-1}-2^v) \le 2^{k-1}-1$; and we have $1\le q\le 2^{h+1}-2^{k-1}$. From here we find the upper and low bonds to $s_i$:
\beas
s_h+0+1\le s_i\le s_{h+1}-s_k+s_{k-1}+2^{h+1}-2^{k-1}\\
s_h <s_i\le s_{h+2}-2^k
\eeas
and due to $k\le h$ we obtain that $s_{h+1}<s_{h+2}-2^k$  and the solution to $s_i$ is fixed: $s_i=s_{h+1}$. Therefore for possible values of  $q$ we get the form $q=s_t+2^v-2^{k-1}, \ t\in\{k,\dots,h\},\ v\in\{0,\dots,k-1\}$. This yields two-parametric series of solutions
\bea\label{disj3}
s_i=s_{h+1}, \quad q=s_t+2^v-2^{k-1}, \ t\in\{k,\dots,h\},\ v\in\{0,\dots,k-1\}.
\eea
It is easy to see that the solutions to $q$ are within segment $[4,s_h]$ (due to $k>2$) and thereby satisfy the restrictions for $q$. \\

The solutions to the forth disjunct\footnote{The 4th disjunct could be included into the 3rd by extending the diapason for $t$ to $[k,h+1]$.} can be described as two-parametric family (with parameters $i,v$). Indeed, 
due to $2^{k-1}-2^v\ge0$ the equation $s_i=(2^{k-1}-2^v)+q$ with restrictions for $v,i,$ and $q$ allows us to list $L$  a series of possible values of $s_i$. First from $s_i=(2^{k-1}-2^v)+q$ we derive $q\le s_i\le q+s_{k-1}$ and with $1\le q\le 2^{h+1}-2^{k-1}$ we can write $1\le s_i\le s_{h+1}$. However in the equation $q=s_i+2^v-2^{k-1}$ the sum $s_i+2^v$ must exceed $2^{k-1}$.  Also when $i=h+1$ due to upper limit for $q$ it should hold that $v=0$. This is why we represent solutions by two series
\bea\label{disj4}
\begin{cases}
q=2^{h+1}-2^{k-1},\ s_i=s_{h+1};&\\
q=s_t+2^v-2^{k-1}, \ s_i=s_t,&  t\in\{k,\dots,h \},\ v\in\{0,\dots,k-1\},
\end{cases}
\eea
where in the second line $i$ and $v$ take their possible values independently. And one more\footnote{Notice that 1 is the minimal value for $v$ here because in case $v=0$ for any value of $i\in\{1,\dots,k-1\}$ the value of $q$ is not positive.}
\beas
 q=s_t+2^v-2^{k-1}, \ s_i=s_t,\  t\in\{1,\dots,k-1\},\ v\in\{1,\dots,k-1\},
\eeas
where $i,v$ are chosen to satisfy the condition $s_i+2^v>2^{k-1}$. The latter condition can be rewritten as $\max\{i,v\}\ge k-1$. Hence we get two sub-series of solutions. One is  $q=s_{k-1}+2^v-2^{k-1}=s_v,\  v=1,\dots k-1,$ whereas the value of  $s_i$ is fixed: $s_i =s_{k-1}$. Another is $q=s_t+2^{k-1}-2^{k-1}=s_t,\ s_i=s_t,\ t=1,\dots,k-1$. The series have a common pair $q=s_i=s_{k-1}$. Therefore we exclude repetition as follows
\bea\label{disj4a}
\begin{cases}
q=s_v,\ s_i=s_{k-1}, & v=1,\dots,k-2;\\
q=s_i,\ s_i=s_t& \ t=1,\dots,k-2;\\
q=s_i=s_{k-1}.&
\end{cases}
\eea
No other solutions exist.\footnote{Notice the similarity between the forms of solutions for the last two disjuncts.}\\

Now we show that for all $q$ s.t. $\Acal_{m+s_{h+1}+2^{k-1}+q}\neq\emptyset$ the number of shifts that solve the equation for $\Acal_{m+s_{h+1}+2^{k-1}+q}$ is even.

Comparing (\ref{disj3}) to the second line of (\ref{disj4}) we find out that solutions for $q$ in both are the same. Therefore for $q=s_t+2^v-2^{k-1},\  t=k,\dots,h,\ v=0,\dots,k-1,$ set $\Acal_{m+s_{h+1}+2^{k-1}+q}$ includes two shifts $s_{t},s_{h+1}$. 

Two first lines from (\ref{disj4a}) imply $\Acal_{m+s_{h+1}+2^{k-1}+s_v}=\{s_v,s_{k-1}\}, v=1,\dots,k-2$.

The third line from (\ref{disj4a}) with the solution $q= s_{k-1},\ s_i=s_{h+1}$ for the second disjunct imply 
$\Acal_{m+s_{h+1}+2^{k-1}+s_{k-1}}=\{s_{k-1},s_{h+1}\}$.

Remaining are the solutions to the first disjunct: $q=s_{h+1}-s_{k-1},\ s_i=s_{h+2}$ and the solution $q=2^{h+1}-2^{k-1},\ s_i=s_{h+1}$ from the 1st line of (\ref{disj4}). Hence $\Acal_{m+s_{h+1}+2_{h+1}}=\{s_{h+1},s_{h+2}\}$. Since $h+1\le n-1$ by the assumption of disjunction on $j$ it is true that $s_{h+1}+2_{h+1}=s_{h+2}\in\Eal_n$. 

Thus we finished the induction on $j$ and the prefix $K_{m,2^n-1}\cdot K_{\eta,\eta}\cdot \Pi_{j=k}^{n-1 } 
 \left[W_k\cdot K_{1,2^j-2^{k-1}+1}\right]$ of $\Eal_n^{2^{n+1}-1}\{K_{m,2^n-1}\}$ is proved.  \\

Now we prove the next block $W_k\cdot K_{1,2^n-2^{k-1}}$. There is one only difference from the previous blocks
$W_k\cdot K_{1,2^j-2^{k-1}+1}$. It is the length of the block $K_{1,2^n-2^{k-1}}$ which is shorter by 1 than in case it was obtained from $W_k\cdot K_{1,2^j-2^{k-1}+1}$ by setting $j=n$. This is enforced by the fact that when $h+1=n$ the previous consideration fails at the very last position $m+s_{h+1}+2_{h+1}$ of zero because when $h+1=n$ the shift $s_{h+2}$ does not belong to $\Eal_n$ and therefore $\Acal_{m+s_{h+1}+2_{h+1}}$ consists of one element $s_{h+1}$ only. 
Simultaneously we see that the position following the block $W_k\cdot K_{1,2^n-2^{k-1}}$ is occupied by 1. It remains only to prove that zeros stay on the remaining $2^n-m-1$ positions. \\

The list of positions of units preceding the suffix $\left[K_{1,2^n-m}\right]\big|_2^{2^n-m}$ is\footnote{Notice the absence of term $2^0$ in the last line. Also $s_{n+1}=s_{n}+2^n$.} 
\beas 
\begin{cases}
m,m+s_k,\\
m+s_k+2^0,\dots, m+s_k+2^u,\dots,  m+s_k+2^{k-1},\\
\dots\\
m+s_{t}+2^0,\dots,m+s_t+2^u,\dots,  m+s_{t}+2^{k-1},\\
\dots\\
m+s_{n}+2^0,\dots,m+s_n+2^u,\dots,  m+s_{n}+2^{k-1},\\
m+s_{n+1}. 
\end{cases}
\eeas
As above we compile an equation for $\Acal_{m+s_{n+1}+q}$ where $1\le q\le 2^n-m-1$:
\beas
s_i\in \Acal_{m+s_{n+1}+q}\iff \hspace{7.5cm}\\
\begin{cases}\notag
 (s_i=s_{n+1}+q) \vee(s_i=s_{n+1}+q-s_k)\ \vee\\
\exists t\in\{k,\dots,n\}\exists v\in\{0,\dots,k-1\}[s_i=s_{n+1}+q-s_t-2^v]\ \vee\\
[s_i=q].
\end{cases}
\eeas
As we know, $2^n-m\le s_k$. Therefore $q\le 2^n-m-1<s_k$. 

The first two disjuncts have no solution for $s_i$ to be in $\Eal_n$ because of $q>0,s_n>s_k, s_{n+1}=s_n+2^n>2\cdot s_n$.

For the 3rd condition we notice that if  $v=0$ then $s_{n+1}-s_t+q-2^v=s_n+2^n-(2^t-1-q+2^0)=s_n+2^n-(2^t-q)$. And since $q>0, t\le n$ we have $2^t-q<2^n$. Therefore $s_i=s_n +2^n-(2^t-q)>s_n$ causing $s_i\notin \Eal_n$. 
Rewrite the third condition as $s_i=(s_{n+1}-2^t)+(q-s_v)$. One series of solutions can be obtained for $t=n$: 
\beas
q=s_v, \ v=1,\dots,k-1, \ s_i=s_n.
\eeas
Let us show that no other solutions exist. Assume $t<n$. In this case $s_{n+1}-s_t+q-2^v> s_n$ because $k\le n, t\le n-1$ and $s_t+2^v-q<2^{n-1}+2^{k-1}-q<2^n$ but $s_{n+1}=s_n+2^n$. Thus we can assume $t=n,\ v\in\{1,\dots,k-1\}$ and $s_i=s_n+(q-s_v)$. Clearly $q>s_v$ implies $s_i\notin \Eal_n$. On the other hand $k\le n, q>0$ entails $q-s_v>-2^{n-1}$. This cancels opportunity for $s_i$ to be lesser $s_n$. 

The 4th disjunct has solutions $q=s_t,\  s_i=s_t,\ t=1,\dots, k-1,$ where the upper limit follows from the condition $q<s_k$.

Thus $\Acal_{m+s_{n+1}+q}$ consists of two shifts $s_t$ and $s_n$ when $q=s_t,\ t\in\{1,\dots, k-1\},$ otherwise is empty. 

\bx

\newpage

\setcounter{theo}{5}

\section*{Addendum 2: a proof of the theorem 6}\label{mult}

\begin{theo}\label{t4}
Let $C=(x_1,\dots,x_r)$ is a position collection, $B$ is its spectrum kernel, and $a\in\Z^+$. Position collection  $aC=(ax_1,\dots,ax_r)$ has  spectrum kernel $\{b\eta(a,b)|b\in B\}$.
\end{theo}
\pr
If $a=1$ nothing is to prove. 
Further, if the spectrum kernel for $C$ is $\{1\}$, the theorem states that the same is true for $aC$ since $\eta(a,1)=1$. And that statement is true because any period $h^k,h\in\{0,1\},$  for $C$ obviously generates a period $h^{ka}$ for $aC$.\footnote{Recall that in contrast to our usage of term ``period'', when we say ``proper period'' we mean a period $\pi$ s.t. no other period $\mu$ exists s.t. $\pi=\mu^m,m>1$.} And yet, both of these periods evidence for the proper period $h$.

Let $a>1, 1\notin B$. Assume $F$ is a semi-infinite word s.t. $F=F\mid_0^{\infty}$ satisfying the recursion set by $aC$. Introduce classes of equivalence  $K_i, i=0,\dots,a-1,$ by modulo $a$ as $K_i=\{x\in\N|x\equiv i\mm a\}$.\footnote{Remember: $\N= \{0\}\cup\Z^+$.} Then, semi-infinite word $F_i$ defined as $F_i(n)=F(an+i)$ satisfies recursion $C$ for each $i$. And vice versa, from any $a$ periods (not necessarily different) of $C$-recursion $\pi_0,\dots,\pi_{a-1}$ it is easy to construct a period of $aC$-recursion. For that on each class $K_i$ we construct a periodic word s.t. $F(an+i)=\pi_i(\rem(n,l_i))$ where $l_i$ is the length of $\pi_i$. Clearly the length of the period of $F$ does not exceed $a\cdot\lcm(l_0,\dots,l_{a-1})$. 
This means that for each $b\in B$ it is possible to construct a period of size $ab$ for the collection $aC$ of shifts. Although it can be not a proper period.

We call $i$-th {\sl component} of $F\mid_0^{ab}$ string $\mathfrak F_i=\Pi_{n=0}^{b-1} F(an+i)$  where $i$
runs over set $\{0,\dots,a-1\}$. Also we write $F=\mathfrak F_0*\frak F_1*\dots*\mathfrak F_{a-1}$.

Now we show how to build for $aC$ a period whose length is $\eta(a,b)b$. Notice that any cyclic permutation of a period is also a period of the same length. 

{\sl Case $a=2$ and $\pi$ is a $C$-period of length $b,b>1,$ where $b\in B$ is an odd integer}:\\
we have $\eta(2,b)=1$ and $\eta(2,b)b=b$, and therefore though 
\beas
F\mid_0^{2b}=\pi(0)\cdot\pi(0)\dots\pi(b-1)\cdot\pi(b-1), 
\eeas
we must show that a cyclic shift of $\frak F_1\ (=\Pi_{n=0}^{b-1}F\mid_0^{2b}(2n+1))$ turns $F\mid_0^{2b}$ into a square of a string $w$ of length $b$. Form here it follows that $b$ is the length of a period for the collection $2C$. 

Since $ww$ should be equal to $\frak F_0*\sigma^x(\frak F_1)$ it must hold $x=\lfloor \frac b2\rfloor $ because just in this case $[\frak F_0*\sigma^x(\frak F_1)](b)=[\frak F_0*\sigma^x(\frak F_1)](0)$. Hence we have expression 
\beas
\pi(0)\cdot\pi(x+1)\cdot\pi(1)\cdot\pi(x+2)\dots\pi(b-1)\cdot\pi(x)\cdot\pi(0)\cdot\pi(x+1)\cdot\pi(1) 
\dots\pi(b-1)\cdot\pi(x) 
\eeas
for $ww$, and it remains to define $w=\pi(0)\cdot\pi(x+1)\cdot\pi(1)\cdot\pi(x+2)\dots\pi(b-1)\cdot\pi(x)$. 

To illustrate the idea of the construction above with more details we consider  
{\sl more general case: $a,b>1, \gcd(a,b)=1$.} 
From $\gcd(a,b)=1$ it follows $\eta(a,b)=1$. We must show that it is possible to choose integers $x_1,\dots,x_{a-1}$ so that 
\beas
\frak F_0*\sigma^{x_1}(\frak F_1)*\dots*\sigma^{x_{a-1}}(\frak F_{a-1})=w^a
\eeas
for some string $w,|w|=b$. To do that we use $a-1$ conditions
\bea\label{conds2}
[\frak F_0*\sigma^{x_1}(\frak F_1)*\dots*\sigma^{x_{a-1}}(\frak F_{a-1})](ib)=[\frak F_0](0), i=\overline{1,a-1}.
\eea
Let $k(i)$ is defined by condition that the position $ib$ belongs to $\frak F_{k(i)}$ i.e. $k(i)=\rem(bi,a)$. Assume that for some $i,i'$ we have $k(i)=k(i')$ or in other words $ib\equiv i'b \mm a$. Since $\gcd(a,b)=1$ that means that $a\mid(i-i')$. The latter contradicts to $0\le i,i'<a$. Thus $k$ is a bijection on $\{1,\dots,a-1\}$  and conditions~(\ref{conds2}) are compatible.  Now we define $x_i=\lfloor\frac{\kappa(i)b}a\rfloor, i=\overline{0,b-1},$ where $\kappa=k^{-1}$.\footnote{By the definition of $F\mid_0^{ab-1}=\frak F_0*\frak F_1*\dots*\frak F_{a-1}$ we have $F(k(i))= \pi(\lfloor\frac{ib}a\rfloor)$. Therefore (to switch on intuition: assume that all letters in $\pi$ are different) $x_{k(i)}= \lfloor\frac{ib}a\rfloor,i\in\{0,\dots,a-1\}$. And yet we need to get rid of this assumption.} 

Let us prove that for chosen $x_1,\dots x_{a-1}$ and any $i\in\{0,\dots,a-1\}$ we have 
\beas
[\frak F_0*\sigma^{x_{1}}(\frak F_1)*\dots*\sigma^{x_{a-1}}(\frak F_{a-1})]\mid_{ib}^{(i+1)b-1}=[\frak F_0*\sigma^{x_1}(\frak F_1)*\dots*\sigma^{x_{a-1}}(\frak F_{a-1})]_0^{b-1}.
\eeas
For that we calculate $[\frak F_0*\sigma^{x_1}(\frak F_1)*\dots*\sigma^{x_{a-1}}(\frak F_{a-1})]_{ib}^{(i+1)b-1}(m),ib\le m<(i+1)b$. 

The number of $a$-block in  $F\mid_0^{ab-1}$ including the position $ib+m$ is  $d=\lfloor \frac{ib+m}a\rfloor$. This means that  $F\mid_0^{ab-1}~(ib+m)=\pi(d)$. On the other hand, position $ib+m$ belongs to component $\frak F_r$ where $r=\rem(ib+m,a)$. However in $\frak F_0*\sigma^{x_1}(\frak F_1)*\dots*\sigma^{x_{a-1}}(\frak F_{a-1})$ components $\frak F_q$ are cyclically shifted to right on $x_q$ positions. And yet, shift $x_r$ is defined by  $j= \kappa(r)$.\footnote{Because the shift $x_r$ is defined by a position $jb$ which belongs to $\frak F_r$.} Since $jb$ and $ib+m$ belongs to the same component\footnote{$\frak F_r$}  the relation $jb\equiv ib+m\mm a$ holds. In other word $\exists s[jb=ib+m+sa]$. 

Thus, the position $d$ in $\pi$ becomes shifted left on $\lfloor\frac{ib+m+sa}a\rfloor$ positions and the actual symbol staying on position $ib+m$ in $[\frak F_0*\sigma^{x_1}(\frak F_1)*\dots*\sigma^{x_{a-1}}(\frak F_{a-1})]_{ib}^{(i+1)b-1}(m)$ is $\pi\left(d-\lfloor\frac{ib+m+sa}a\rfloor\right)= \pi\left(\lfloor \frac{ib+m}a\rfloor- \lfloor\frac{ib+m+sa}a\rfloor\right)=\pi(\rem(-s,b))$.

It remains to show that $\rem(-s,b)$ does not depend on $i$. But we have (see above) $jb=ib+m+sa$ which means that $b\mid(m+sa)$. Since $j$ is defined uniquely there should exist one only $s$ for each $m\in\{0,\dots,b-1\}$. However this formally does not cancel opportunity that $s$ depends on $i,j$ too. On the other hand, if $s,s'$ obeys $b\mid(m+sa)$ and $b\mid(m+s'a)$ then $b\mid(s-s')$. And $s,s'<b$ holds because $sa,s'a<ab$ due to the fact that we operate with the segment of $F$ of length $ab$ exactly. Hence we can find $s$ only on the basis of $a,b,m$ does not matter what value has $i$. \\

{\sl The general case}. We can assume $a,b, \eta(a,b)=\eta>1$.  
Let 
\bea\label{F}
F=\pi(0)\dots\pi(0)\pi(1)\dots\pi(1)\dots \pi(b-1)\dots\pi(b-1)
\eea
where each group of symbols $\pi(i)\dots\pi(i)$ has length $a$.

Let us denote by $\mu$ the integer $\frac a{\eta(a,b)}$. We must show that it is possible to choose integers $x_1,\dots,x_{a-1}$ so that 
\bea\label{hatF}
\hat F=\frak F_0*\sigma^{x_1}(\frak F_1)*\dots*\sigma^{x_{a-1}}(\frak F_{a-1})=w^{\mu}
\eea
for some string $w,|w|=\eta b$. To do that we use some conditions from $ab$ of them
\bea\label{conds3}
\hat F(j)=\hat F(i\eta b+j),0= 1,\dots,\mu-1,j=0,\dots\eta b-1
\eea
and then check the other. First we choose $a$ equations
\bea\label{conds4}
\hat F(0)=\hat F(i\eta b), i= 0,\dots,\mu-1
\eea
and define mapping
\beas
k:\{0,1,\dots,\mu-1\}\to\{0,\eta\dots,(\mu-1)\eta\}  
\eeas
as 
\bea\label{ki}
k(i)=\rem(i\eta b,a). 
\eea
Indeed, $\eta\mid k(i)$ for any $i\in \{0,1,\dots,\mu-1\}$ because $\eta\mid a$. Also $k$ is a bijection. To see the latter assume that $k(i)=k(i')$ for some $i,i'\in\{0,\dots, \mu-1\}$. From (\ref{ki}) we deduce $i\eta b\equiv i'\eta b\mm a$ or $a\mid (i-i')\eta b$. Since $a=\eta\mu$ and $\gcd(\mu,b)=1$ it follows from here that $\mu\mid(i-i')$. However $|i-i'|<\mu$. Hence $i=i'$. It remains to notice that $|\{0,1,\dots,\mu-1\}|=|\{0,\eta\dots,(\mu-1)\eta\}|$.

Now we define 
\bea\label{x_k(i)}
x_{k(i)}=\left\lfloor\frac{i\eta(a,b)b}{a}\right\rfloor, i=0,\dots,\mu-1.
\eea 
In this way we have defined $\mu$ shifts from $a$ of them  $x_0,\dots,x_{a-1}$ in (\ref{hatF}). 
To define others, for each $j\in\{0,\dots,a-1\}$ we set 
\bea\label{x_j}
x_j=x_{\bar j}, \text{ where } j\in\{0,\dots,a-1\}, \bar j=\left\lfloor\frac j{\eta}\right\rfloor\eta=j-\rem(j,\eta).
\eea
The definition is correct because $\bar j\in\{0,\eta,\dots,(\mu-1)\eta\}$ and thereby $\bar j$ is a value of function $k$; so we can apply (\ref{x_k(i)}) using $k^{-1}$ as 
\bea\label{x_j}
x_j=\left\lfloor\frac{k^{-1}(\bar j)\eta b}{a}\right\rfloor, j=0,\dots,a-1.
\eea

Thus $\hat F$ is completely defined and we should show that $\hat F(i\eta b+j),0\le j<\eta b,$ does not depend on $i$.
The position $i\eta b+j$ belongs to component $\frak F_{r}$ where $r = \rem(i\eta b+j,a)$. Therefore $F(i\eta b+j)= \pi(\lfloor\frac{i\eta b+j}a\rfloor)$ (see (\ref{F})). 

Let us for beginning assume that $j=\bar j$, i.e $\eta\mid j$. In this case $r\in\{0,\eta,\dots,\eta(\mu-1)\}$, $r=\bar r$,  $k^{-1}(\bar r)$ is defined and equal to $k^{-1}(r)$.

In $\hat F$ the component $\frak F_r$ is cyclically shifted right by $x_r$ positions and therefore $\hat F(i\eta b+j)=\pi(\lfloor\frac{i\eta b+j}a\rfloor-x_r)$ where $x_r= \left\lfloor \frac{k^{-1}(r)\eta b}{a}\right\rfloor$.
For some $i'$ we can write $k^{-1}(r)\eta b=i'\eta b$. 
An integer $g=g(j,b,a)$ exists s.t. $k^{-1}(i\eta b+j)\equiv(i+g)\eta b\mm a$. Indeed, within equivalence by modulo $a$ the amount can be rewritten as $i\eta b+g\eta b$ where $j\equiv g\eta b\mm a$. When $j=\eta j'$ the latter equation can be rewritten as $j'\equiv gb\mm{\mu}$ which has a solution $g$ because $j$ could be first replaced with $\rem(j,a)$ which yields for $j'$ the bounds $0\le j'<\mu$ and because $\gcd(b,\mu)=1$.\footnote{Since $\gcd(b,\mu)=1$ there are integers $h,q$ s.t. $1=qb+h\mu$. From here $j'=j'qb+j'h\mu$. From here when $0\le j'<\mu$ we deduce $j'=j'qb\mm{\mu}$. So $g=\rem(qj',\mu)$.\\  Example: $a=15,b=3,\eta=3,\mu=5$. In this case we have:
\beas
j=1\cdot\eta\implies k^{-1}(i\eta b+j)\equiv(i+2)\eta b\mm a\\
j=2\cdot\eta\implies k^{-1}(i\eta b+j)\equiv(i+4)\eta b\mm a.
\eeas
}
On the other hand, positions $i\eta b+j, k^{-1}(r)\eta b$ belong to the same component $\frak F_r$. Therefore  $i\eta b+j- (i+g)\eta b=sa$ for some integer $s$. In other words we have $sa=j-g\eta b$. The latter entails 
\beas
\left\lfloor\frac{i\eta b+j}a\right\rfloor- \left\lfloor \frac{k^{-1}(r)\eta b}{a}\right\rfloor=\frac{j-g\eta b}a.
\eeas
This means that the shift does not depend on $i$.\footnote{In the example (see the previous footnote) we had $g=2$ for $j=\eta$ and $g=4$ for $j=2\eta$. Respectively the fraction $\frac{j-g\eta b}a$ is equal to $\frac {3-2\cdot 9}{15}=-1$ and 
$\frac {6-4\cdot 9}{15}=-2$ }

Further we need to expand the proof on the case when $\eta\nmid j$. It is possible to partition $\hat F$ into sub-strings of length $\eta$ because $a=\eta\mu$. By the same reason the partition consists of partitions of all $a$-blocks. The latter means that $j$ belongs to an $a$-block iff the whole $\eta$-block $(\bar j,\bar j+1,\dots,\bar j+\eta-1)$ belongs to the $a$-block. By definition (\ref{x_j}) and setting $F$ we have  $x_j=x_{j'}$ if $\bar j=\bar{j'}$ and $F(j)=F(j')$. The latter implies $\hat F(j)=\hat F(j')$. So, if $\hat F(i\eta b+\bar j)$ does not depend on $i$, $\hat F(i\eta b+j)$ does not depend on $i$ as well.

It remains to prove that $aC$ has no period not divisible by some $b'\eta(a,b'),b'\in B$.
Assume there exists a minimal period of length $q$ for  collection $aC$ in the sense that $q$ is not divisible by the length of more short period for $aC$. Then there exists a period of length $\frac{\lcm(q,a)}a$ for collection $C$. Since $\lcm(q,a)\le qa$ we have $\frac{\lcm(q,a)}a\le q$. And since any period length for collection $C$ is divisible by some $b\in B$ from $\frac{\lcm(q,a)}a\mid q$ it follows that $b\mid q$.

We need to deduce that $\frac q{\gcd(q,a)}$ is a period for $C$. (Here $\gcd(q,a)$ plays a role of $\eta(a,b)$.) 
For that show that $\frac a{\gcd(a,q)}$ and $\gcd(a,q)$ relatively prime.

When $\gcd(a,q)>1$, either $\frac a{\gcd(a,q)}$  or $\frac q{\gcd(a,q)}$ are relatively prime with $\gcd(a,q)$ but only one of them. If the first then $\gcd(a,q)$ is in fact $\eta(a,b)$ and $q=\eta(a,b)b$ where $b=\frac q{\gcd(a,q)}$. If the second, then this $b$ is relatively prime with $a$ and as it was shown above $b$ is also a period length for $aC$. Hence, $q$ is not a minimal (also in the sense that is not divisible by other period lengths). 

If $\gcd(a,q)=1$ then $q$ is also a period length for $C$.

\bx

\ED